\documentclass[11pt, a4paper]{article}
\usepackage{amsmath}
\usepackage{amsfonts}
\usepackage{amssymb}
\usepackage{ulem}

\usepackage{amsthm}
\usepackage{url}
\usepackage{dcolumn}

\usepackage{cite}
\usepackage{fancyhdr}

\usepackage{graphicx, color}
\usepackage{hyperref}
\usepackage{xcolor}

\newtheorem{thm}{Theorem}

\newtheorem*{rem*}{Remark}

\usepackage{here}

\begin{document}

\title{Convex pentagonal monotiles in the 15 Type families}
\author{ Teruhisa SUGIMOTO$^{ 1), 2)}$ }
\date{}
\maketitle

{\footnotesize

\begin{center}
$^{1)}$ The Interdisciplinary Institute of Science, Technology and Art

$^{2)}$ Japan Tessellation Design Association

E-mail: ismsugi@gmail.com
\end{center}

}

{\small
\begin{abstract}
\noindent
The properties of convex pentagonal monotiles in the 15 Type families and their tilings are 
summarized. The Venn diagrams of the 15 Type families are also shown.
\end{abstract}
}


\section{Introduction}
\label{section1}

If all tiles in a tiling\footnote{ 
A \textit{tiling} (or \textit{tessellation}) of the plane is a collection of sets, called tiles, 
that cover the plane without gaps or overlaps, except for the boundaries of the tiles.
The term ``tile" refers to a topological disk, whose boundary is a simple closed 
curve \cite{G_and_S_1987}.
} 
are of the same size and shape, the tiling is described as \textit{monohedral}, which 
allows the use of reflected tiles (posterior side tiles) in monohedral tilings \cite{G_and_S_1987, 
wiki_pentagon_tiling}. In other words, in monohedral tilings, the anterior and posterior sides of 
the tiles are treated as the same type (i.e., the concept assumes that there is only one 
type of tile). In this study, a tile that admits a monohedral tiling is referred to as a \textit{monotile}.

\begin{figure}[hb]
\centering\includegraphics[width=14.5cm,clip]{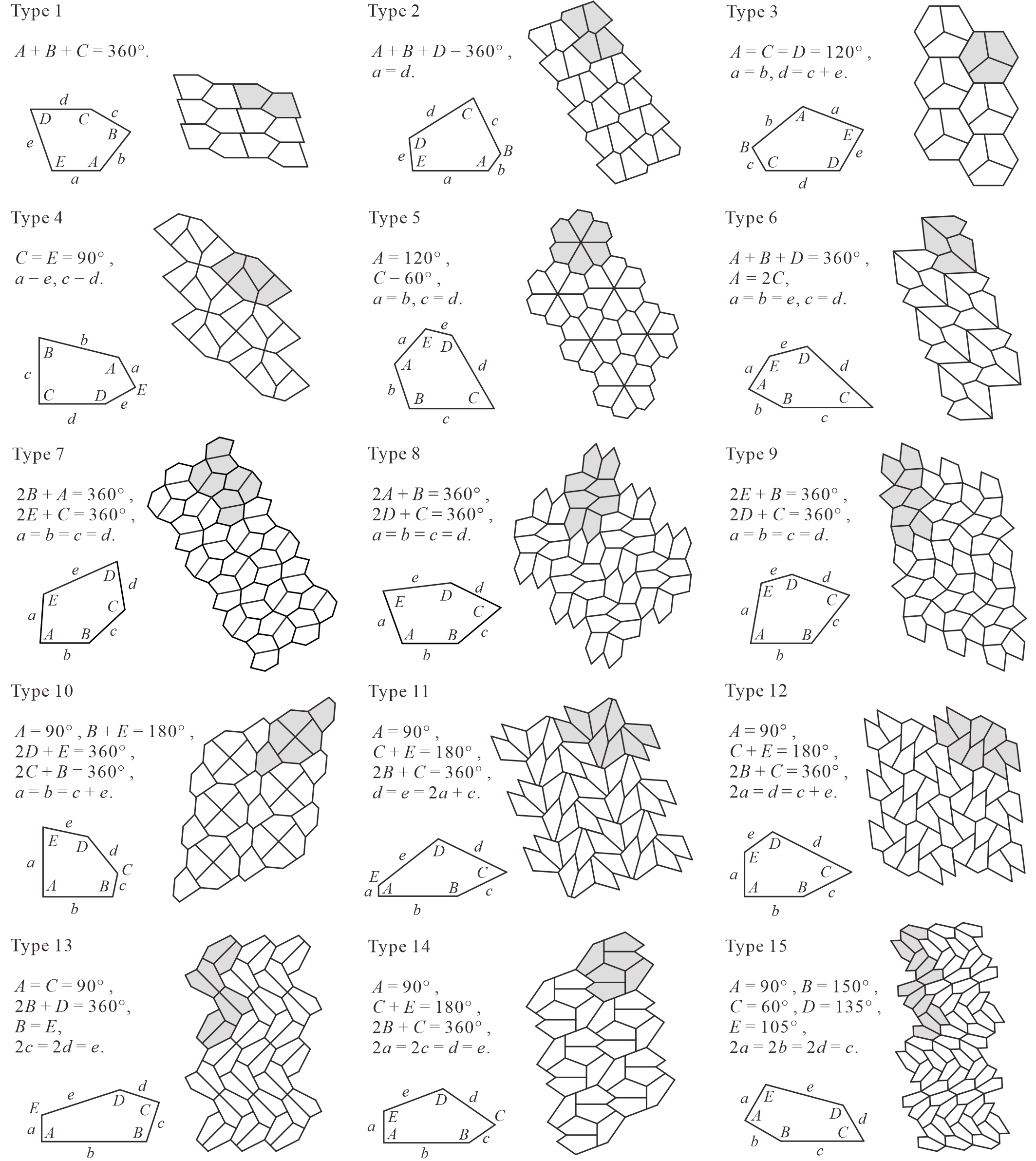} 
\caption{{\small 
Fifteen Type families of convex pentagonal monotiles.
}
\label{fig01}
}
\end{figure}

Currently, there are 15 known families of convex pentagonal monotiles, each labeled as a ``Type,'' 
as shown in the list of Figure~\ref{fig01} \cite{Gardner_1975a, Gardner_1975b, G_and_S_1987, Kershner_1968, 
Mann_2015, Rao_2017, Reinhardt_1918, Schattschneider_1978, Scherphuis, Stein_1985, 
Sugimoto_2012, Sugimoto_2015a, Sugimoto_2015b, Sugimoto_2016, Sugimoto_2017a, Sugimoto_2018a, 
Sugimoto_2017b, Sugimoto_2018b, wiki_pentagon_tiling, Zong_2020}. Each convex pentagonal 
monotile in the list is defined by the conditions expressing the fundamental relationship between 
the edge lengths and angles of the tile; however, some degrees of freedom remain. For example, 
convex pentagonal monotiles belonging to the Type 1 family (also referred to as convex 
pentagonal monotiles belonging to Type 1, or simply as convex pentagonal monotiles of 
Type 1) satisfy the property that the sum of three consecutive angles is $360^ \circ $ 
or that the sum of the remaining two consecutive angles is $180^ \circ $. In Figure~\ref{fig01}, 
this property (relationship) is expressed as $A + B + C = 360^ \circ $, 
which represents the tile conditions of Type 1 (the Type 1 family) (see \ref{appA}). 
Each convex pentagonal monotile  belonging to the Types 14 and 15 families has no degrees of 
freedom other than size (that is, each convex pentagonal monotile of Types 14 and 15 has 
a fixed shape. For example, the value of angle $C$ in the convex pentagon of Type 14 is 
$\cos ^{ - 1} \bigl( \bigl( 3\sqrt {57} - 17 \bigr) / 16 \bigr) \approx 1.2099\;\mbox{rad} \approx 69.32^ \circ$). 
Because some degrees of freedom other than size may remain, these Type families of 
convex pentagonal monotiles are not necessarily ``disjoint.'' The tilings depicted for each 
convex pentagon shown in Figure~\ref{fig01} are representative of each Type (Type family). 
A representative tiling can be formed using only the relationships derived from the 
tile conditions of each Type. The gray region in each tiling in Figure~\ref{fig01} indicates 
a \textit{translation unit} (a unit that can generate periodic\footnote{ 
A tiling exhibits \textit{periodicity} if its translation by a non-zero vector coincides with 
itself; a tiling is considered periodic if it coincides with its translation by two 
linearly independent vectors. However, in this study, a tiling with periodicity is 
referred to as \textit{periodic}, and a tiling without periodicity is referred to as 
\textit{non-periodic} \cite{Sugimoto_2024}.
} 
tiling through translation alone)\footnote{ 
Typically, a translation unit of a tiling is chosen to consist of the minimum number of tiles, 
as illustrated by each tiling in Figure~\ref{fig01}. 
}. 
Thus, the convex pentagonal monotiles in these 15 Type 
families can generate periodic tilings, as the representative tiling of each Type 
includes a translation unit, as shown in Figure~\ref{fig01}.

The purpose of this study is to summarize the properties of the convex pentagonal 
monotiles and tilings of the 15 Type families depicted in Figure~\ref{fig01}. 
In particular, we show the ``disjoint/not disjoint" relationship of each Type family 
of convex pentagonal monotiles, and specifically illustrate convex pentagonal monotiles 
belonging to multiple Type families in cases where they are not disjoint.

\renewcommand{\figurename}{{\small Figure}}


\section{Properties of convex pentagonal monotiles and their tilings}
\label{section2}

Let $T_{x}$ denote a set of convex pentagonal monotiles belonging to a Type $x$ 
family for $x = 1 - 15$, as shown in Figure~\ref{fig01}.

\begin{figure}[t]
\centering\includegraphics[width=15cm,clip]{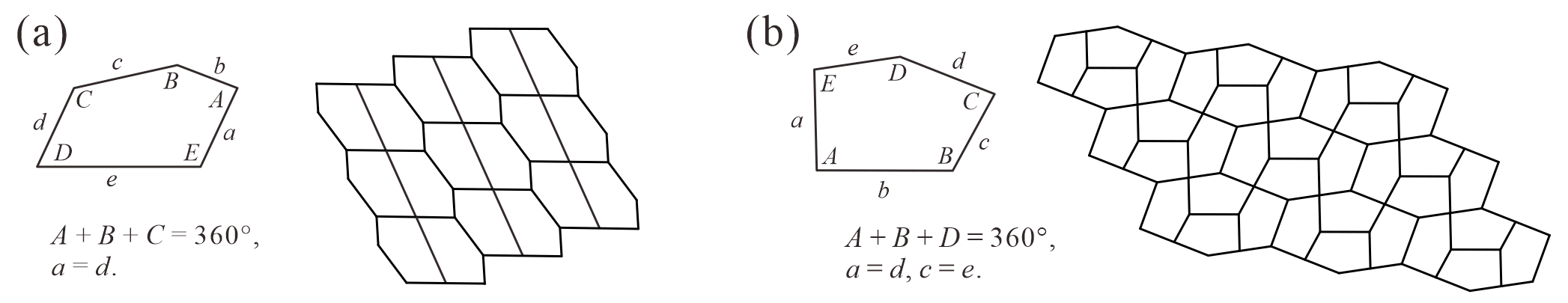} 
\caption{{\small 
Examples of edge-to-edge tilings with convex pentagonal monotiles 
belonging to the Type 1 or Type 2 families.
}
\label{fig02}
}
\end{figure}

A tiling with convex polygons is defined as \textit{edge-to-edge} if any two convex polygons 
in the tiling are either disjoint or share one vertex or an entire edge of the tiling in 
common\footnote{ 
An edge-to-edge tiling with polygons is defined as a tiling in which the vertices (corners) and 
edges (sides) of the polygons coincide with the vertices (points where three or more tiles meet) 
and edges of the tiling \cite{G_and_S_1987, Schattschneider_1978}.
}. 
Otherwise, the tiling is classified as \textit{non-edge-to-edge} \cite{G_and_S_1987, Sugimoto_2012, 
Sugimoto_2015a, Sugimoto_2015b, Sugimoto_2016, Sugimoto_2017a, Sugimoto_2017b}. 
While the representative tilings of Type 1 or Type 2 shown in Figure~\ref{fig01} are 
non-edge-to-edge, $T_{1}$ or $T_{2}$ contain convex pentagonal monotiles that can 
generate edge-to-edge tilings. For example, as shown in Figure~\ref{fig02}(a), the convex 
pentagonal monotile satisfying the conditions ``$A + B + C = 360^ \circ ,\;a = d$" 
belongs to the Type 1 family and can form an edge-to-edge tiling. (Remark: this convex 
pentagonal monotile can also generate non-edge-to-edge tilings.) Similarly, as shown in 
Figure~\ref{fig02}(b), the convex pentagonal monotile satisfying the conditions 
``$A + B + D = 360^ \circ ,\;a = d, \;c = e$" belongs to the Type 2 family and can form 
an edge-to-edge tiling \cite{Sugimoto_2012, Sugimoto_2015b, Sugimoto_2017a, Sugimoto_2017b}. 
By contrast, for instance, the convex pentagonal monotiles in $T_{3}$ or $T_{13}$ cannot 
generate edge-to-edge tilings.

\renewcommand{\figurename}{{\small Figure}}

In \cite{Sugimoto_2016}, we presented the following Theorem regarding convex pentagonal 
monotiles that can generate edge-to-edge tilings\footnote{ 
We have known that the same result as Theorem~\ref{thm1} was obtained by Bagina 
\cite{Bagina_2011, Bagina_2012} in 2011 after we derived Theorem~\ref{thm1} in 2012.
}:

\begin{thm}\label{thm1}
If a convex pentagon can generate an edge-to-edge monohedral tiling, it belongs to at 
least one of the Type 1, 2, or 4--9 families shown in Figure~\ref{fig01}.
\rm{(}\textit{In other words, convex pentagonal monotiles that can generate an 
edge-to-edge tiling are contained in $T_{1}$, $T_{2}$, $T_{4}$, $T_{5}$, $T_{6}$, $T_{7}$, 
$T_{8}$, or $T_{9}$.}\rm{)}
\end{thm}

As mentioned above, the Type families of convex pentagonal monotiles are not necessarily 
disjoint. In other words, some convex pentagonal monotiles belong to multiple Type families 
(i.e., some convex pentagonal monotiles are contained in multiple $T_{x}$).

\renewcommand{\figurename}{{\small Figure}}
\begin{figure}[htbp]
\centering\includegraphics[width=9cm,clip]{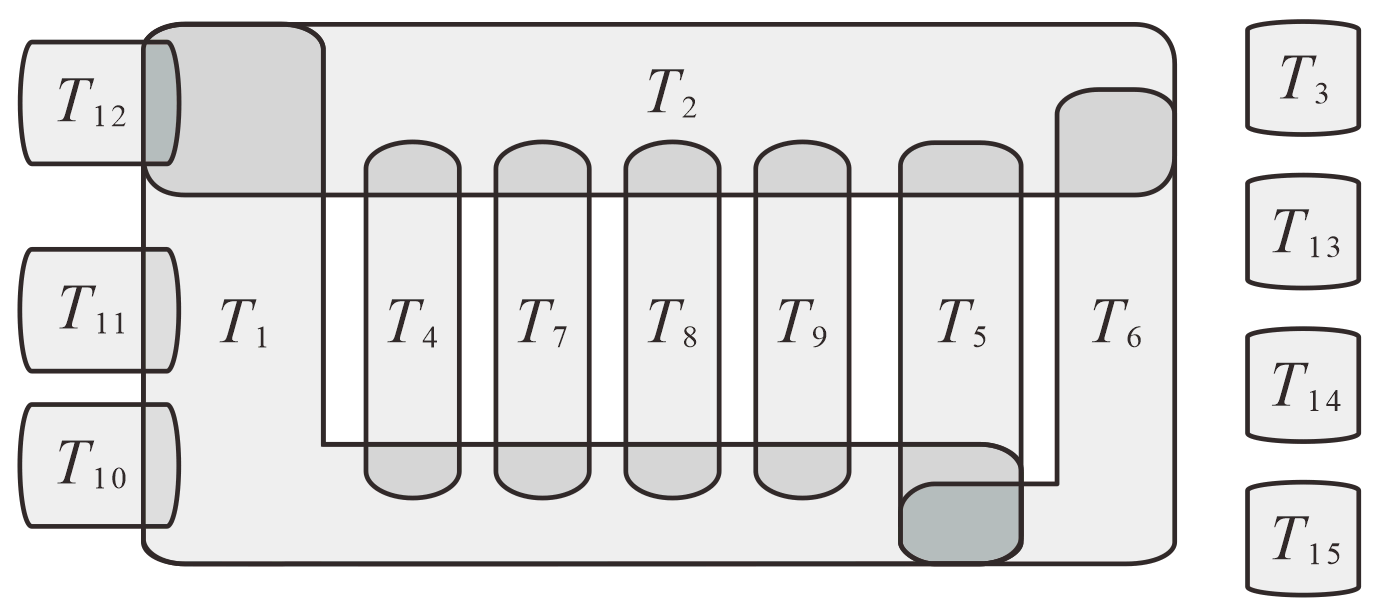} 
\caption{{\small 
Venn diagram of $T_{x}$.
}
\label{fig03}
}
\end{figure}

\renewcommand{\figurename}{{\small Figure}}
\begin{figure}[htbp]
\centering\includegraphics[width=15cm,clip]{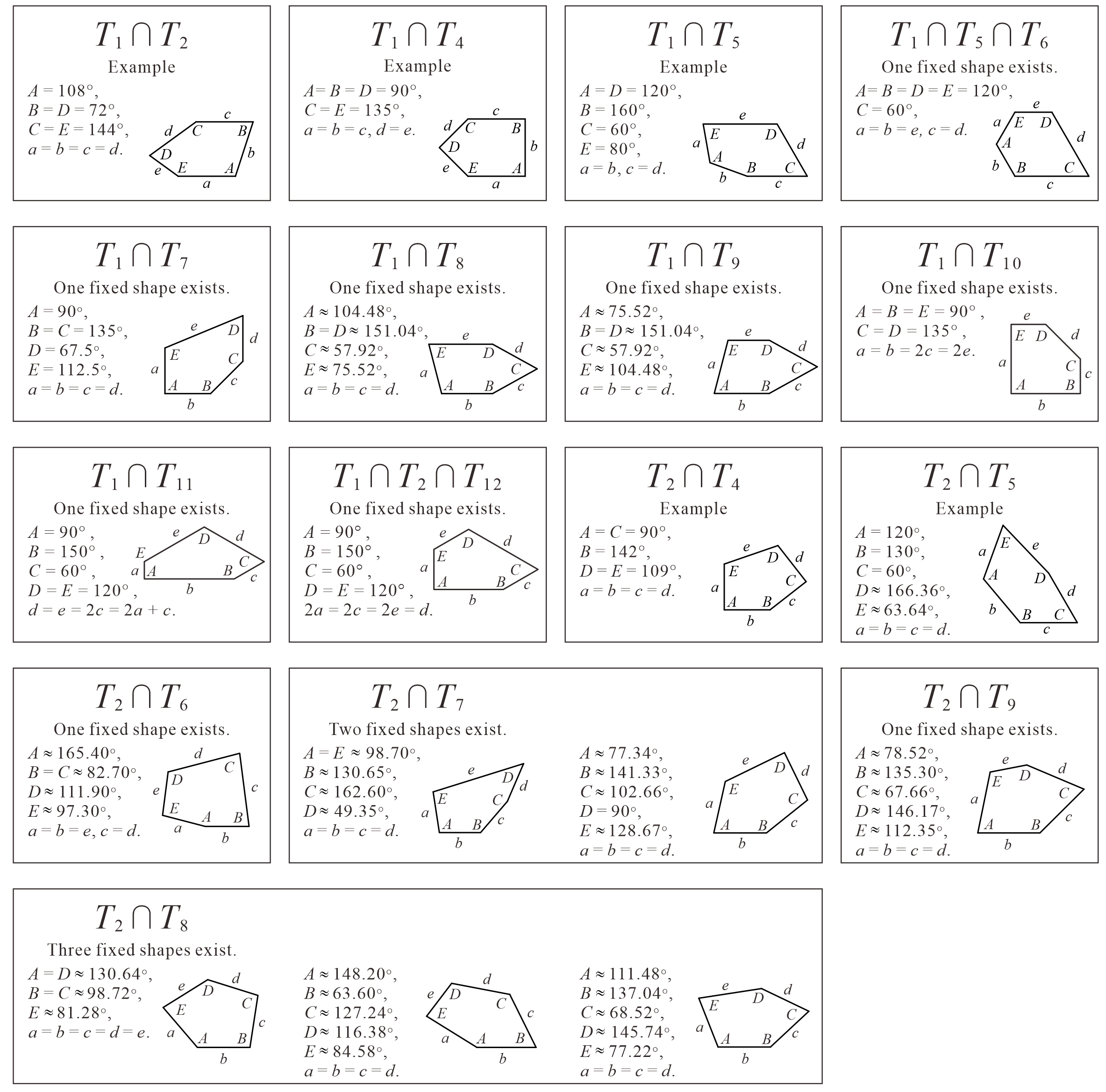} 
\caption{{\small 
Convex pentagonal monotiles contained in each intersection of the 
Venn diagram shown  in Figure~\ref{fig03}.
}
\label{fig04}
}
\end{figure}

Figure~\ref{fig03} shows the Venn diagram of $T_{x}$. Figure~\ref{fig04} presents the convex 
pentagonal tiles that are contained in each intersection of the Venn diagram in 
Figure~\ref{fig03} (see \ref{appA}). Various shapes exist for each convex pentagonal 
monotile in the intersections ``$T_{1} \cap T_{2}$, $T_{1} \cap T_{4}$, $T_{1} \cap T_{5}$, 
$ T_{2} \cap T_{4}$, and $T_{2} \cap T_{5}$." The convex pentagonal monotiles of their 
intersections shown in Figure~\ref{fig04} are examples of each case (hence the label 
``Example" in the figure). In contrast, each of the intersections ``$T_{1} \cap T_{5} \cap T_{6}$, 
$T_{1} \cap T_{7}$, $T_{1} \cap T_{8}$, $ T_{1} \cap T_{9}$, $T_{1} \cap T_{10}$, 
$T_{1} \cap T_{11}$, $T_{1} \cap T_{2} \cap T_{12}$, $T_{2} \cap T_{6}$, and $T_{2} \cap T_{9}$" 
contains one fixed shape. The intersection ``$T_{2} \cap T_{7}$" contains two fixed shapes, 
and the intersection ``$T_{2} \cap T_{8}$" contains three fixed shapes. Convex pentagonal 
monotiles belonging to the Type 3 or Type 13 families have degrees of freedom other than 
size. Therefore, the convex pentagonal monotiles in $T_{3}$ or $T_{13}$ exhibit various 
shapes; however, these convex pentagons do not belong to any other Type families. In contrast, 
$T_{14}$ or $T_{15}$ each contain one fixed shape, and these convex pentagons also 
do not belong to any other Type families.

\begin{figure}[tb]
\centering\includegraphics[width=15cm,clip]{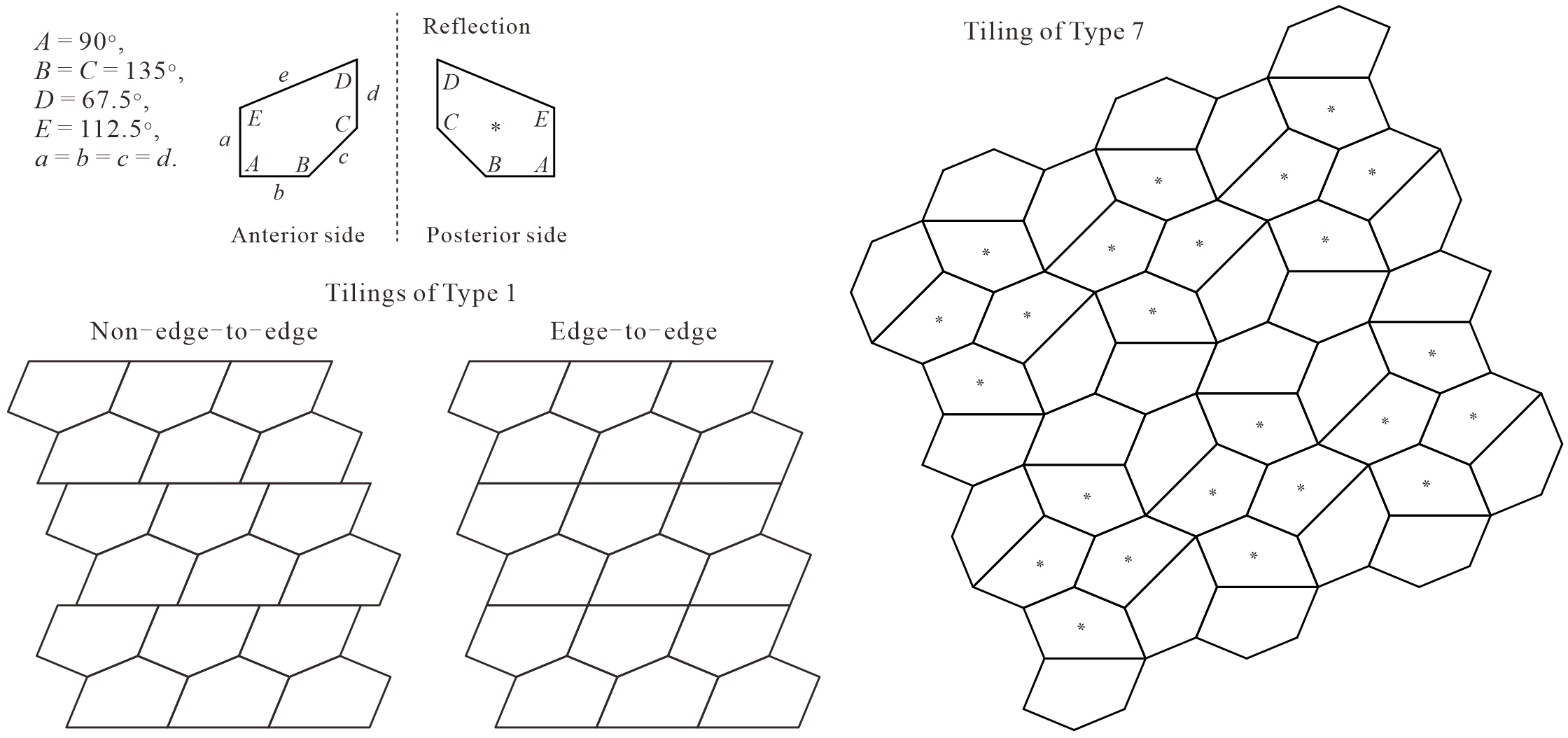} 
\caption{{\small 
Convex pentagonal monotile belonging to both the Type 1 and Type 7 
families, the representative tilings of Type 1 and of Type 7 with the convex 
pentagon. In this figure, the convex pentagons corresponding to the 
reflected tiles are marked with an asterisk ``*."
}
\label{fig05}
}
\end{figure}

For example, because the convex pentagonal monotile in $T_{1} \cap T_{7}$ belongs to 
both the Type 1 and Type 7 families, it can form representative tilings of Types 1 and 7, 
as shown in Figure~\ref{fig05} \cite{Sugimoto_2020d}.

\renewcommand{\figurename}{{\small Figure}}

\begin{figure}[tb]
\centering\includegraphics[width=15cm,clip]{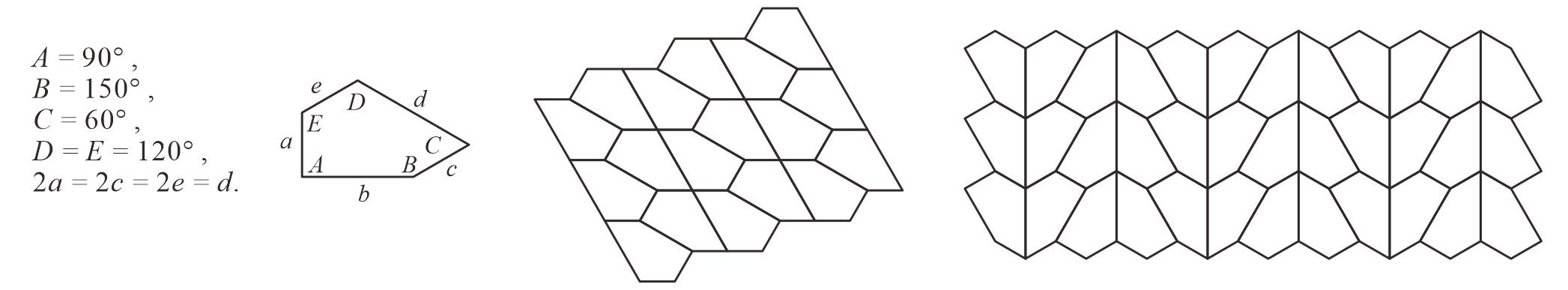} 
\caption{{\small 
Examples of edge-to-edge tilings formed by convex pentagonal 
monotile in $T_{1} \cap T_{2} \cap T_{12}$.
}
\label{fig06}
}
\end{figure}

The convex pentagonal monotile in $T_{1} \cap T_{2} \cap T_{12}$ (i.e., the convex pentagonal 
monotile belonging to the Type 1, Type 2, and Type 12 families) can form edge-to-edge 
tilings as shown Figure~\ref{fig06}. In other words, the set of convex pentagonal 
monotiles that can generate edge-to-edge tilings contains a convex pentagonal 
monotile belonging to the Type 12 family. (Among the convex pentagonal monotiles 
contained in $T_{12}$, only the convex pentagon that also belongs to the Type 1 and 
Type 2 families can generate edge-to-edge tilings.) This does not negate Theorem~\ref{thm1}.

As shown in Figure~\ref{fig01} (or Figures~\ref{fig13} and \ref{fig14} in \ref{appD}), 
the representative tilings of Types 1 and 3--6 
do not use reflected tiles (convex pentagons). In other words, if a convex 
pentagonal monotile belongs to at least one of the Type 1 or 3--6 families 
shown in Figure~\ref{fig01}, it can generate tilings without the use of reflected 
tiles (that is, convex pentagonal monotiles in $T_{1}$, $T_{3}$, $T_{4}$, $T_{5}$, 
or $T_{6}$ can generate tilings without the use of reflected tiles).  By contrast, 
as shown in Figure~\ref{fig01}, the representative tilings of Types 2 and 7--15 use 
reflected tiles. For example, if $A = 86^ \circ $ for the tile conditions of Type 7 in 
Figure~\ref{fig01}, the convex pentagonal monotile with $A = 86^ \circ $ 
belongs only to the Type 7 family and can only form the representative 
tiling of Type 7, as depicted in Figure~\ref{fig01}. Only the relationships 
``$2B + A = 360^ \circ ,\; 2E + C = 360^ \circ ,\; 2D + A + C = 360^ \circ $" 
are used for the vertices of this tiling, and reflected tiles are required to 
form this tiling. This means that there are monotiles that cannot generate a tiling if 
the use of reflected tiles is not allowed during the tiling generation process, even if 
they are monotiles (see \ref{appB}). (This classification highlights that some monotiles 
require reflections for successful tiling generation, while others do not.) 
It should be noted that among the convex pentagonal monotiles contained in $T_{7}$, 
the convex pentagon that also belongs to the Type 1 family (i.e., the convex pentagonal 
monotile belonging to both the Type 1 and Type 7 families) can generate tilings without 
the use of reflected tiles, as shown in Figure~\ref{fig05}. That is, the tilings of Type 1 
depicted in Figure~\ref{fig05} do not use reflected tiles.

\renewcommand{\figurename}{{\small Figure}}

\begin{figure}[ht]
\centering\includegraphics[width=14.5cm,clip]{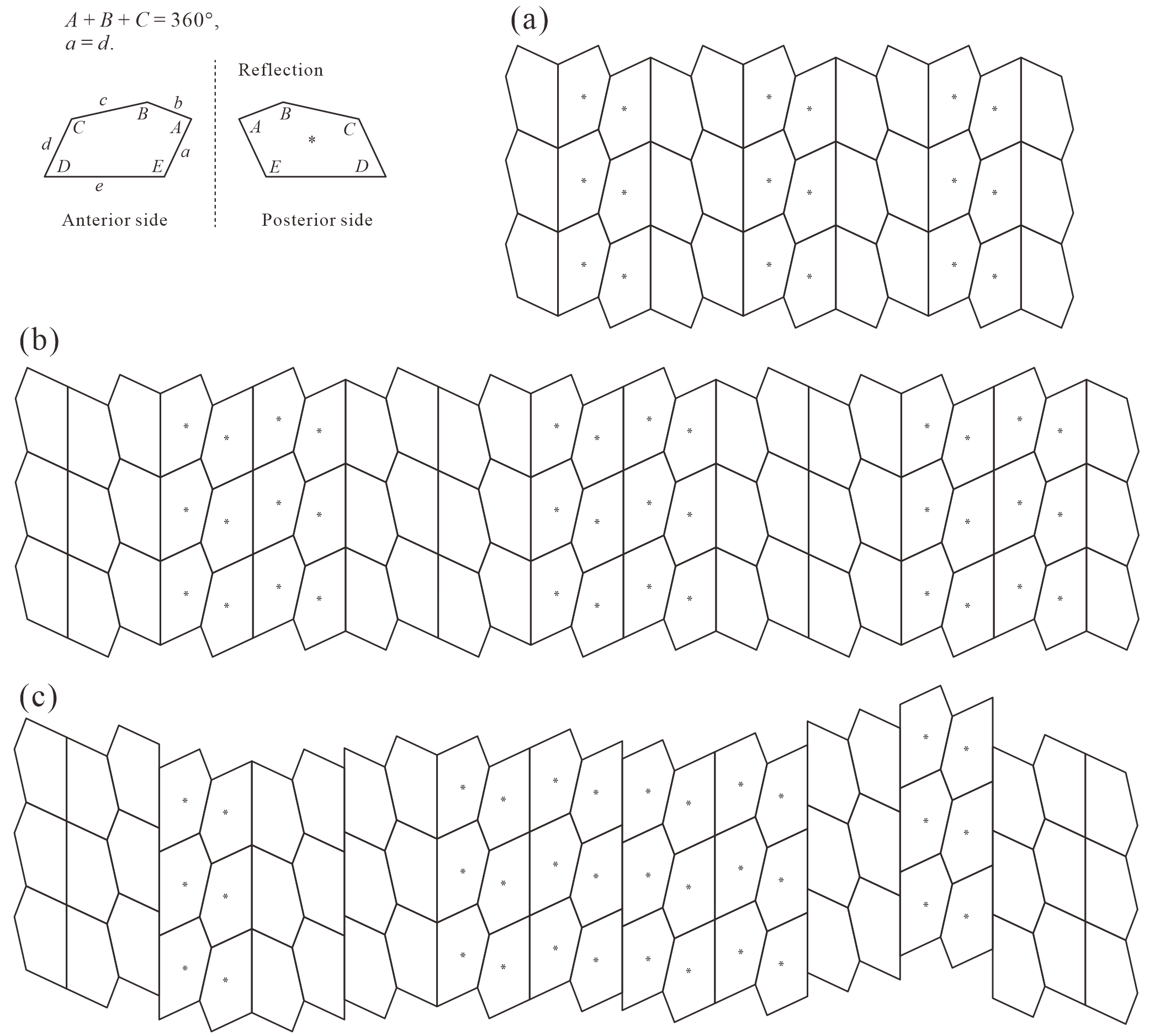} 
\caption{{\small 
Tilings that are formed by convex pentagonal monotile belonging to the 
Type 1 family. In this figure, the convex pentagons corresponding to the 
reflected tiles are marked with an asterisk ``*."
}
\label{fig07}
}
\end{figure}

\begin{figure}[t]
\centering\includegraphics[width=14.5cm,clip]{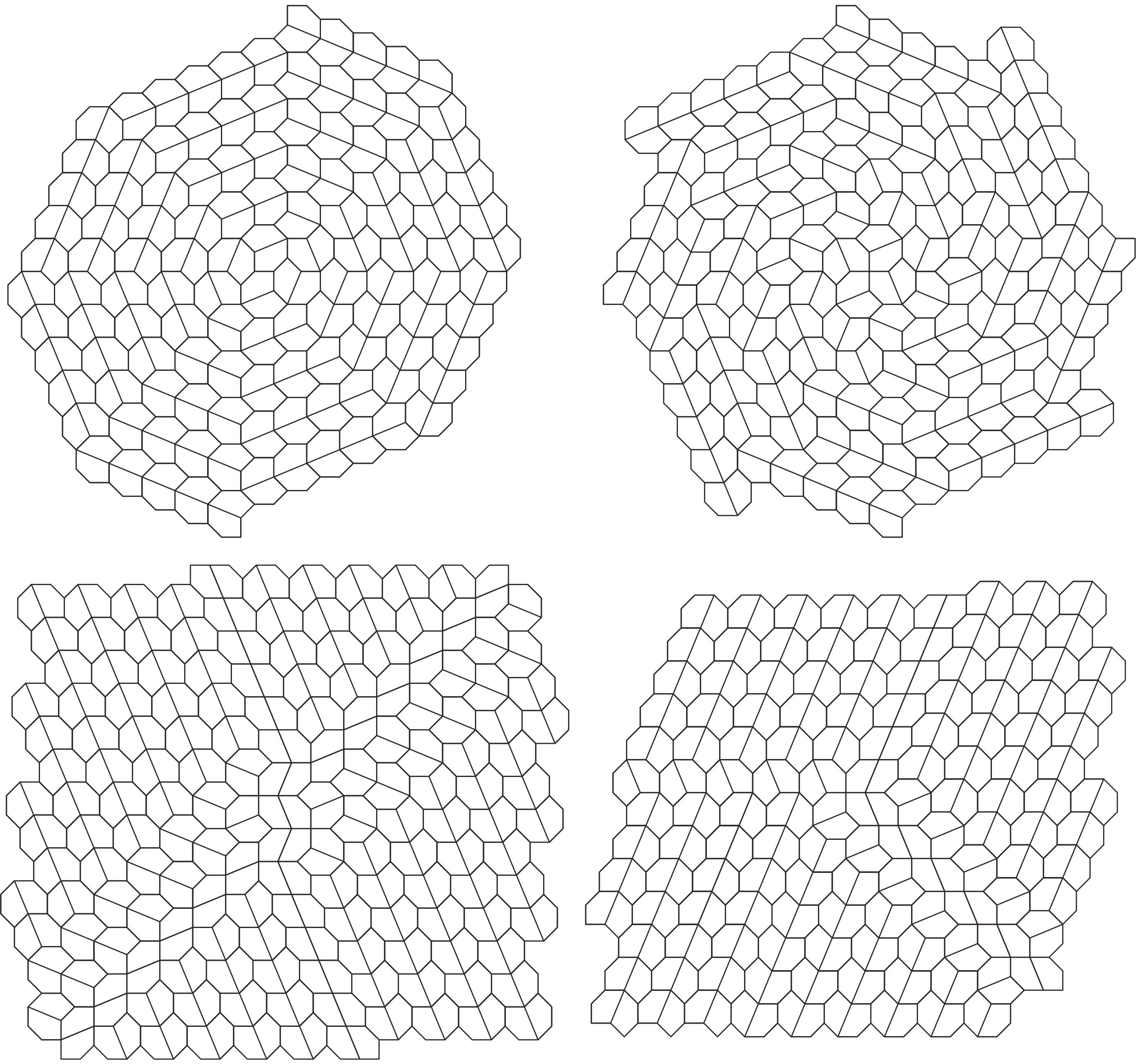} 
\caption{{\small 
Non-periodic tilings formed by convex pentagonal monotile 
belonging to both the Type 1 and Type 7 families.
}
\label{fig08}
}
\end{figure}

Some convex pentagonal monotiles can form tiling patterns (design patterns created 
by the arrangement of polygonal tiles) that differ from the representative tiling of each 
Type \cite{Scherphuis, Sugimoto_2020a, Sugimoto_2020b, Sugimoto_2020c, Sugimoto_2020d, 
Sugi_Araki_2017a, Sugi_Araki_2017b}. For example, the convex pentagonal monotile in 
Figure~\ref{fig02}(a), which belongs to the Type 1 family and satisfies ``$A + B + C = 360^ \circ ,
\; a = d$," can form tilings such as those illustrated in Figures~\ref{fig07}(a) and \ref{fig07}(b) 
by using reflected tiles. Additionally, the convex pentagon satisfying 
``$A + B + C = 360^ \circ ,\;a = d$" can generate tilings, as shown in Figure~\ref{fig07}(c), 
using reflected tiles because the belts it forms through vertical translation in the same 
direction can be freely connected horizontally\footnote{ 
In such tilings, various translation units can be formed, and it is also possible to generate 
tilings with no translation units.
}. 
Figure~\ref{fig08} presents examples of tilings formed by the convex pentagonal 
monotile that belongs to both the Type 1 and Type 7 families \cite{Sugimoto_2020d}. 
As shown in Figure~\ref{fig08}, some convex pentagonal monotiles can generate 
non-periodic tilings \cite{Sugimoto_2020a, Sugimoto_2020b, Sugimoto_2020c, 
Sugimoto_2020d, Sugi_Araki_2017a, Sugi_Araki_2017b}. (Remark: As noted earlier, 
the convex pentagonal monotiles in the 15 Type families can always generate periodic tilings).

Some tiling patterns cannot be formed without the use of reflected tiles. 
However, the tiles forming such tiling patterns do not necessarily need to 
be reflected to generate the tilings.

\renewcommand{\figurename}{{\small Figure}}

\renewcommand{\figurename}{{\small Figure}}

Let $\Im $ denote a tiling with congruent tiles on the Euclidean plane. 
A \textit{symmetry} of $\Im $ is an isometry\footnote{ 
Isometries include rotation, translation, reflection, glide reflection, and the identity 
transformation \cite{G_and_S_1987, Weisstein_isometry}.
} 
of the Euclidean plane that maps the tiles of $\Im $ onto tiles of $\Im $. The 
\textit{symmetry group} of $\Im $, denoted as $\mbox{Sym}(\Im )$, is the collection 
of all such symmetries. Two tiles $P_{1}$ and $P_{2}$ of $\Im $ are said to be 
\textit{equivalent} if there exists a symmetry belonging to $\mbox{Sym}(\Im )$ that 
maps $P_{1}$ onto $P_{2}$. The collection of all tiles of $\Im $ that are 
equivalent to $P_{1}$ is called the \textit{transitivity class} of $P_{1}$. If $\Im $ 
has precisely $k$ transitivity classes, then it is said to be \textit{$k$-isohedral} or 
\textit{tile-$k$-transitive}. Isohedral means 1-isohedral. 
In addition, a monotile that cannot generate isohedral tilings is referred to 
as an \textit{anisohedral tile} (i.e., an anisohedral tile is a prototile\footnote{ 
See \ref{appC} for information on ``prototile."
}
 that admits monohedral tilings but no isohedral tilings) \cite{Goodman_2017, 
G_and_S_1987, Mann_2015, Scherphuis, Weisstein_anisohedral, Weisstein_isohedral, Zong_2020}.

The representative tilings of Types 1--5 shown in Figure~\ref{fig01} are isohedral. 
In contrast, the representative tilings of Types 6--9 and 11--13 shown in Figure~\ref{fig01} 
are 2-isohedral, while those of Types 10, 14, and 15 are 3-isohedral (see \ref{appD}) 
\cite{G_and_S_1987, Mann_2015, Schattschneider_1978, Scherphuis, wiki_pentagon_tiling, 
Zong_2020}. 
 Convex pentagonal monotiles in $T_{1}$, $T_{2}$, $T_{3}$, $T_{4}$, or $T_{5}$ (i.e., convex 
 pentagonal monotiles belonging to at least one of the Type 1--5 families shown in 
 Figure~\ref{fig01}) can generate isohedral tilings \cite{G_and_S_1987}. 
However, note that some convex pentagonal monotiles in $T_{1}$, $T_{2}$, 
$T_{3}$, $T_{4}$, or $T_{5}$ can also generate non-isohedral tilings, such as 
2-isohedral tilings. Consequently, convex pentagonal monotiles that are not 
contained in $T_{1}$, $T_{2}$, $T_{3}$, $T_{4}$, or $T_{5}$ are classified as 
anisohedral tiles. (It should be noted that some convex pentagonal monotiles 
in $T_{7}$, $T_{8}$, $T_{9}$, $T_{10}$, $T_{11}$, or $T_{12}$ are not anisohedral 
tiles. As shown in Figure~\ref{fig04}, the convex pentagonal monotiles in the 
intersections ``$T_{1} \cap T_{5} \cap T_{6}$, $T_{1} \cap T_{7}$, 
$T_{1} \cap T_{8}$, $ T_{1} \cap T_{9}$, $T_{1} \cap T_{10}$, 
$T_{1} \cap T_{11}$, $T_{1} \cap T_{2} \cap T_{12}$, $T_{2} \cap T_{6}$, 
$T_{2} \cap T_{7}$, $T_{2} \cap T_{8}$, and $T_{2} \cap T_{9}$'' are 
not anisohedral tiles.) From the Venn diagram in Figure~\ref{fig03}, all convex 
pentagonal monotiles in $T_{13}$, $T_{14}$, or $T_{15}$ are classified as 
anisohedral tiles.


\section{Conclusion}
\label{section3}

In their work\cite{Mann_2015}, Mann et al. propose the following conjecture: 
``It would be reasonable to conjecture that any 
unmarked convex pentagon that admits a tiling of the plane admits at least one periodic 
tiling; that is, it would be reasonable to conjecture that there are no aperiodic\footnote{ 
See \ref{appC} for information on ``aperiodic."
} 
convex pentagons. If this conjecture is true, then all convex pentagons that admit 
tilings of the plane also admit at least one \textit{i}-block transitive tiling. Thus, the class 
of pentagons being studied in this article may well encompass all possible pentagons 
that admit tilings of the plane." 

Rao \cite{Rao_2017} claimed that all convex pentagonal monotiles belong to the known 15 
Type families. Therefore, if his claim is correct, convex pentagonal monotiles can always 
generate periodic tilings. However, it appears that Rao's claim has not yet been 
definitively established \cite{Rao_2017, wiki_pentagon_tiling}.

Based on the studies of Mann et al. and Rao, we also consider that convex pentagonal 
monotiles that can generate periodic tilings belong to the known 15 Type families.


\medskip
\noindent
\textbf{Acknowledgments.} 
The authors would like to thank Professor Emeritus Yoshio Agaoka, 
Hiroshima University, for providing valuable comments.
\bigskip



\bigskip

\appendix
\def\thesection{Appendix \Alph{section}}
\section{}
\label{appA}

This section supplements the concept of tile conditions for the Type families of convex 
pentagonal monotiles and their tilings.

\begin{figure}[htbp]
\centering\includegraphics[width=15cm,clip]{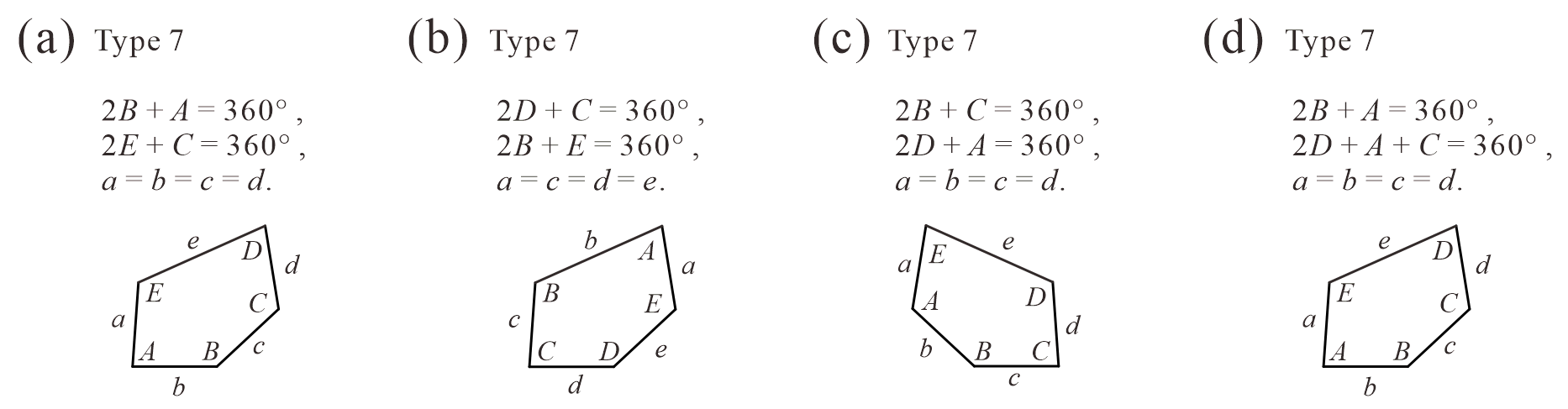} 
\caption{{\small 
Examples of various notations for the tile conditions of convex 
pentagonal monotiles belonging to the Type 7 family.
}
\label{fig09}
}
\end{figure}

Figure~\ref{fig09} illustrates the tile conditions (i.e., the conditions expressing the fundamental 
relationship between the edge lengths and angles) of a convex pentagonal monotile belonging 
to the Type 7 family using four different notations. Figures~\ref{fig09}(a) (with the same notation 
as in Figure~\ref{fig01}) and \ref{fig09}(b) show the same tile conditions, but with the symbols 
``$A,\;B,\;C,\;D,\;E$" for the vertices (angles) and ``$a,\;b,\;c,\;d,\;e$" for the edges of the 
convex pentagons positioned differently. The convex pentagon in Figure~\ref{fig09}(c) is a 
reflection image of the convex pentagon in Figure~\ref{fig09}(a), and it represents the same 
tile conditions according to the vertex and edge symbols assigned to it. Additionally, 
although one of the angle relationships differs between Figures~\ref{fig09}(a) and 
\ref{fig09}(d), the relationships are mathematically equivalent because the sum of the 
interior angles of a convex pentagon equals $540^ \circ $. Thus, these 
relationships can be derived from one another.

\renewcommand{\figurename}{{\small Figure}}

\renewcommand{\figurename}{{\small Figure}}
\begin{figure}[h]
\centering\includegraphics[width=15cm,clip]{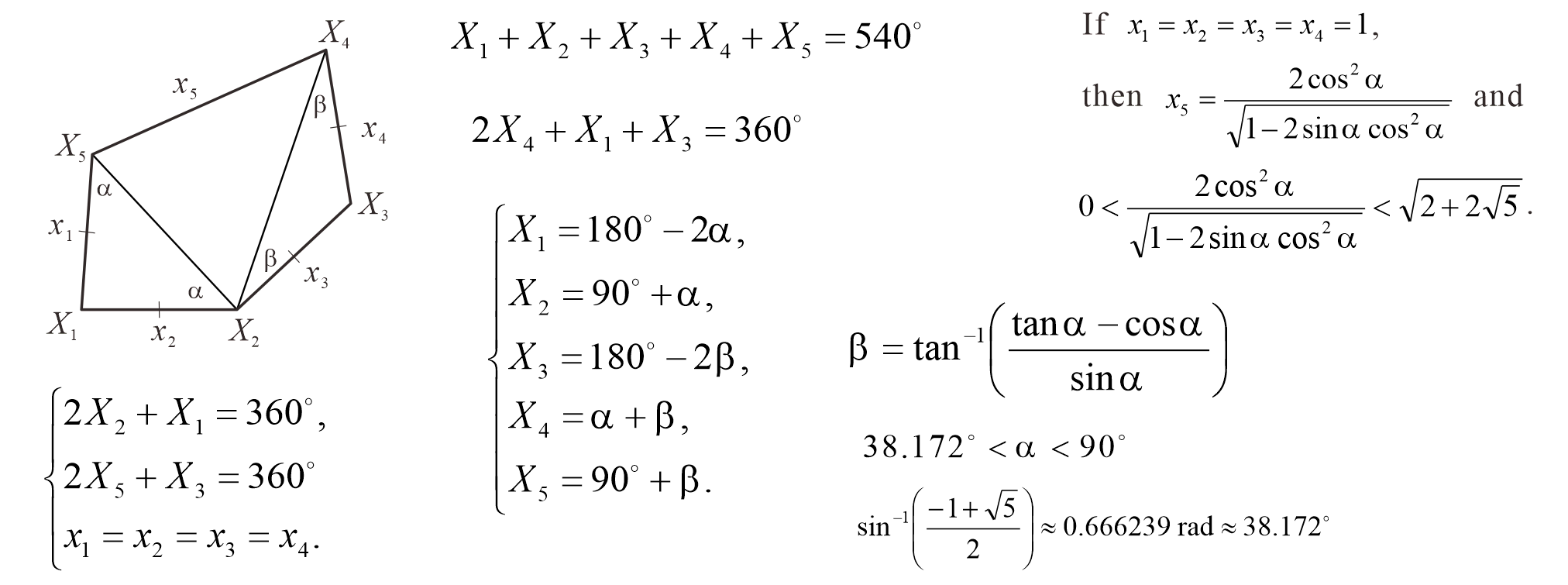} 
\caption{{\small 
Geometric properties of convex pentagonal monotiles belonging to 
the Type 7 family.
}
\label{fig10}
}
\end{figure}

In conclusion, convex pentagonal monotiles belonging to the Type 7 family are defined as 
convex pentagons with the geometric properties shown in Figure~\ref{fig10}. (As indicated in 
Figure~\ref{fig10}, if a convex pentagon with vertex labels ``$X_1 ,\;X_2 ,\;X_3 ,\;X_4 ,\;X_5 $" 
and edge labels ``$x_1 ,\;x_2 ,\;x_3 ,\;x_4 ,\;x_5 $" satisfies the relationships 
``$2X_2 + X_1 = 360^ \circ ,\;2X_5 + X_3 = 360^ \circ ,\;x_1 = x_2 = x_3 = x_4 $," the other 
geometric properties in Figure~\ref{fig10} can be derived based on the sum of the interior 
angles of the convex pentagon.) As shown in Figure~\ref{fig01}, we use the notation in 
Figure~\ref{fig09}(a) to indicate convex pentagonal monotiles belonging to the Type 7 family; 
however, we could have used a different notation (e.g., those in Figures~\ref{fig09}(b), 
\ref{fig09}(c), and \ref{fig09}(d)). Therefore, although the tile conditions (i.e., expressions 
of the formulas) may differ depending on the notations in the various lists of convex 
pentagonal monotiles,  it is essential to understand that these notations refer to the 
same tile conditions.

\begin{figure}[t]
\centering\includegraphics[width=5.1cm,clip]{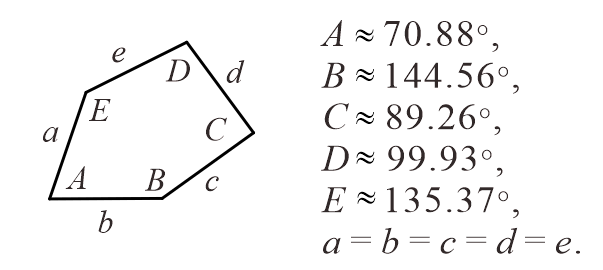} 
\caption{{\small 
Equilateral convex pentagonal monotile belonging to the Type 7 
family.
}
\label{fig11}
}
\end{figure}

The shape of the convex pentagon in Figure~\ref{fig10} is determined by assigning a 
value to $\alpha$, for example. Therefore, convex pentagonal monotiles belonging to the 
Type 7 family are considered to have one degree of freedom other than size. By adjusting 
the remaining degrees of freedom, we identify the equilateral case in the convex pentagonal 
monotiles belonging to the Type 7 family, as shown in Figure~\ref{fig11}. (The notation 
of the tile conditions of Type 7 requires them to have four equal edge lengths but does 
not prohibit them from having five equal edge lengths.) The equilateral convex pentagonal 
monotile of the Type 7 family has no degrees of freedom other than its size; therefore, 
it has a fixed shape. Furthermore, Figure~\ref{fig01} shows that the tile conditions of 
Types 8 and 9 are similar to those of Type 7. Based on the geometric properties, 
we observe an equilateral convex pentagonal monotile belonging to the Type 8 family 
that also belongs to the Type 2 family\footnote{ 
Hirschhorn and Hunt presented a theorem to the effect that ``If an equilateral convex 
pentagon can generate a monohedral tiling, it belongs to at least one of the Type 1, 2, 
or 7 families shown in Figure~\ref{fig01}" \cite{Bagina_2004, Hirschhorn_1985}.
} 
(see the equilateral convex pentagon of the intersection ``$T_{2} \cap T_{8}$" in 
Figure~\ref{fig04}). By contrast, we observe that there is no equilateral convex pentagonal 
monotile belonging to the Type 9 family. (It would appear that the notation of the tile 
conditions of Type 9 requires them to have four equal edge lengths but does not prohibit 
them from having five equal edge lengths. However, owing to the geometric properties, it 
cannot have five equal edge lengths) \cite{Schattschneider_1978}.

\renewcommand{\figurename}{{\small Figure}}

\renewcommand{\figurename}{{\small Figure}}
\begin{figure}[t]
\centering\includegraphics[width=15cm,clip]{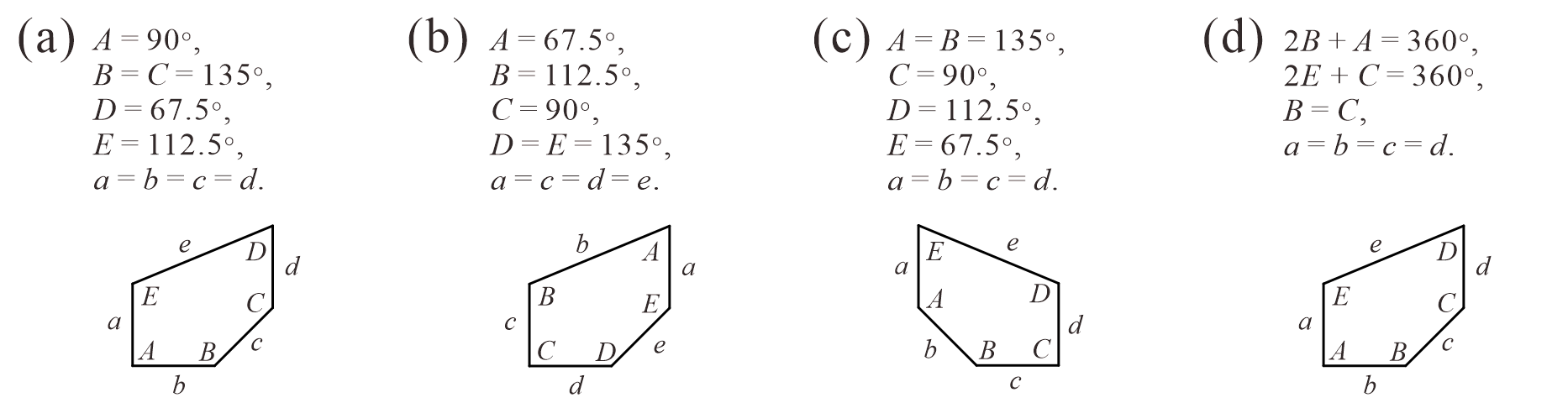} 
\caption{{\small 
Examples of various notations for the tile conditions of convex 
pentagonal monotile belonging to both the Type 1 and Type 7 families.
}
\label{fig12}
}
\end{figure}

Considering the geometric properties of Figure~\ref{fig10}, we find that there exists 
one fixed shape, the convex pentagonal monotile, which belongs to both the 
Type 1 and Type 7 families\footnote{ 
For the convex pentagon in Figure~\ref{fig10}, it can be observed that the property that 
``the sum of three consecutive angles is $360^ \circ $ (= the tile conditions of Type 1)" 
does not hold for all values of $\alpha$.
} 
(see the convex pentagon of the intersection 
``$T_{1} \cap T_{7}$" in Figure~\ref{fig04}). The relationship between the edge lengths 
and angles of the convex pentagon (i.e., the tile conditions) can be expressed in various 
manners, as shown in Figure~\ref{fig12}. The notation in Figure~\ref{fig12}(a) is adjusted 
such that the positions of the edge and vertex symbols correspond directly to the 
notations of the tile conditions of Types 1 and 7 in Figure~\ref{fig01}. In contrast, 
if the values of angles $A$ and $B$ in Figures~\ref{fig12}(b) or \ref{fig12}(c) are 
directly substituted into the angle relationship ``$2B + A = 360^ \circ $" in 
Figure~\ref{fig09}(a), this clearly does not hold. (It is important not to erroneously 
conclude that the convex pentagon in Figure~\ref{fig12} does not belong to the 
Type 7 family simply because the relationships in the formulas in the list do not hold 
with such a simple substitution.) For example, the notation of Figure~\ref{fig12}(b) 
confirms that the tile conditions of Type 7 are satisfied by the notation of 
Figure~\ref{fig09}(b), and the tile conditions of Type 1 are satisfied from 
$C + D + E = 360^ \circ$. Although the notation in Figure~\ref{fig12}(d) does 
not indicate the concrete values of each angle, these values (the same as those 
in Figure~\ref{fig12}(a)) can be derived from the relationship formulas.

The appearance of the formula for the tile conditions expressed in the list 
of convex pentagonal monotiles will change depending on how the vertex and 
edge symbols of the convex pentagon are assigned, how the relationship 
between the vertices and edge lengths is selected, and how the ``anterior 
side/posterior side" of the convex pentagon is set. ``Type" is a classification 
based on the essentially different geometric properties possessed by convex 
pentagons. Therefore, to determine the Type family to which a certain convex 
pentagonal monotile belongs, the geometric properties of the convex pentagon 
need to be considered. 

There are many convex pentagonal monotiles that can form tiling patterns that differ 
from the representative tilings of Type \cite{Scherphuis, Sugimoto_2020a, Sugimoto_2020b, 
Sugimoto_2020c, Sugimoto_2020d, Sugi_Araki_2017a, Sugi_Araki_2017b, wiki_pentagon_tiling}. 
In particular, convex pentagonal monotiles that belong to multiple Type families have 
the potential to form various tiling patterns \cite{Sugimoto_2020c, Sugimoto_2020d, 
Sugi_Araki_2017a, Sugi_Araki_2017b}. Thus, if a convex pentagon forms a tiling 
pattern that differs from the representative tilings in the list of 15 Type families 
of convex pentagonal monotiles, one should not immediately assume 
that it is an unknown convex pentagonal monotile. In the list, we usually adjust 
special cases, such as convex pentagonal monotiles belonging to multiple Type 
families, so that they are not used as example figures of convex pentagonal 
monotiles or representative tilings belonging to each Type family. For example, 
in the case of Type 7 in Figure~\ref{fig01}, we present the tiling whose nodes 
(vertices of the tiling) are formed only by the relationships 
``$2B + A = 360^ \circ ,\;2E + C = 360^ \circ ,\;2D + A + C = 360^ \circ $" 
that can be derived from the tile conditions, and we present the convex pentagon 
whose tiling can be generated only by these relationships. (If a new tiling pattern with 
convex pentagonal monotiles is thought to be found, simply because the 
tiling pattern cannot be formed by the convex pentagons traced from the 
list, it is premature to assume that the convex pentagon forming the tiling 
pattern does not belong to the known Type family.) Furthermore, the classification 
of the Type families and the classification of tiling patterns are different. 
As the tiling patterns that can be formed by convex pentagonal monotiles are 
innumerable, it is difficult to determine whether a tiling pattern is already known.


\section{}
\label{appB}

In this section, we consider convex pentagonal monotiles that require the use 
of reflected tiles during the tiling generation process. However, in this discussion, 
convex pentagons with line symmetry are excluded because they do not have a 
distinction between their anterior and posterior (reflected) sides.

As mentioned in Section~\ref{section2}, convex pentagonal monotiles in $T_{1}$, $T_{3}$, 
$T_{4}$, $T_{5}$, or $T_{6}$ can generate tilings without the use of reflected tiles. 
Considering the Venn diagram in Figure~\ref{fig03} and the list in Figure~\ref{fig04}, 
convex pentagonal monotiles that are not contained in $T_{1}$, $T_{3}$, $T_{4}$, 
$T_{5}$, or $T_{6}$ correspond to convex pentagonal monotiles belonging only 
to the Type 2 family, convex pentagonal monotiles belonging to the Types 
7--12 families that are not contained in $T_{1}$, and convex pentagonal 
monotiles belonging to the Types 13--15 families\footnote{ 
The convex pentagons contained in the intersections ``$T_{2} \cap T_{7}$, 
$T_{2} \cap T_{8}$, and $T_{2} \cap T_{9}$" in Figures~\ref{fig03} and \ref{fig04} 
correspond to convex pentagonal monotiles belonging to the Types 7--12 families 
that are not contained in $T_{1}$. The convex pentagons in the other intersections 
correspond to convex pentagonal monotiles in $T_{1}$, $T_{3}$, $T_{4}$, 
$T_{5}$, or $T_{6}$. Note that the equilateral convex pentagonal monotile in 
$T_{2} \cap T_{8}$ is excluded because this convex pentagon exhibits line symmetry 
(see Figure~\ref{fig04}).
}. 
Among these, because each convex pentagonal monotile in the Types 14 and 15 families 
has a fixed shape and can only generate the specific tilings of each type shown in 
Figure~\ref{fig01}, it is evident that they cannot generate tilings without the use of 
reflected tiles. In contrast, it currently appears unclear whether all convex pentagonal 
monotiles without line symmetry in the remaining cases can generate tilings without 
the use of reflected tiles. Here, let $T_{only - 2}$ denote the set of convex pentagonal 
monotiles belonging only to the Type 2 family, and let $T_{7to13}$ denote the set of 
convex pentagonal monotiles belonging to the Types 7--13 families. Then, the set of 
the above ``remaining cases" can be expressed as 
$T_{only - 2} \cup \bigl( T_{7to13} \cap T_{1}^{C} \bigr)$. 
Note that $T_{x}^{C}$ denotes the complement of $T_{x}$ (the set of convex pentagonal 
monotiles not belonging to the Type $x$ family).

As described in \ref{appA}, the shape of convex pentagonal monotiles belonging to the 
Type 7 family can be modified by adjusting the remaining degrees of freedom. Consequently, 
it can be observed that cases exist in which the sum of certain interior angles equals 
$360^ \circ $, in addition to the angle relationships ``$2B + A = 360^ \circ ,\;2E + C = 360^ \circ 
,\;2D + A + C = 360^ \circ $" that are derived from the tile conditions of Type 7 in 
Figure~\ref{fig01}. For example, the convex pentagonal monotiles in $T_{2} \cap T_{7}$ 
shown in Figure~\ref{fig04} represent such cases. Clearly, these monotiles can form 
the representative tilings of Types 2 and 7, both of which require the use of reflected tiles. 
However, it remains unclear whether they can generate other tiling patterns, and if so, 
whether such patterns would also require the use of reflected tiles. 

If a convex pentagonal monotile in $T_{x}$ exhibits properties in which the 
sum of two or more interior angles equal $360^ \circ $ or $180^ \circ $, in 
addition to those derived from the tile conditions of Type $x$, it may be able 
to generate tiling patterns other than the representative tilings of Type $x $ \footnote{ 
A convex pentagonal monotile in $T_{x}$ can always form the representative 
tilings of Type $x$. However, it is not immediately clear from the fundamental 
relationships between the edge lengths and angles of the tile whether other 
tiling patterns can be generated. Furthermore, no complete classification is 
available of which convex pentagonal monotiles in $T_{x}$ can generate which 
tiling patterns.
}. 

A search is conducted for all convex pentagonal monotiles in 
$T_{only - 2} \cup \bigl(T_{7to13} \cap T_{1}^{C} \bigr)$ that can generate tiling patterns 
other than the representative tilings of Types 2 and 7--13. If all such tiling patterns 
require the use of reflected tiles and the tiles do not exhibit line 
symmetry, the following conclusion can be drawn: ``If convex pentagonal 
monotiles in the 15 Type families do not exhibit line symmetry and do not 
belong to the Type 1 or 3--6 families  shown in Figure~\ref{fig01}, then the convex 
pentagonal monotiles cannot generate tilings without the use of reflected tiles.''

\section{}
\label{appC}

In general, when considering tilings, there are only a finite number of tile shapes 
that can be used. These finite-shape diagrams are called \textit{prototiles}, 
and a set of prototiles admits tilings of the plane. The tiling generated by a set of 
prototiles can cover the plane infinitely. In a monohedral tiling, because all the tiles 
are the same size and shape (i.e., congruent), there is only one type of prototile. 
A set of prototiles is said to be \textit{aperiodic} (i.e., an ``aperiodic set of prototiles") 
if copies of the prototiles can be assembled into tilings of the plane such that all 
tilings with the prototiles are non-periodic \cite{Hallard_1991, G_and_S_1987, 
wiki_aperiodic_set, wiki_list_aperiodic}.

The Penrose tiles are a well-known example of aperiodic set of prototiles. They are 
considered set consisting of two polygons with matching conditions \cite{Hallard_1991, 
G_and_S_1987, Sugimoto_2017a, wiki_aperiodic_set, wiki_list_aperiodic}. 
Matching conditions specify how tiles must connect to form a valid tiling, which 
can sometimes be represented by assigning colors or orientations to specific 
edges of the prototiles.

A key problem regarding aperiodic set of prototiles, often referred to as the ``Einstein problem'', 
asks:  ``Is there a single aperiodic prototile (with or without a matching condition), 
that is, one that admits only non-periodic tilings by congruent copies?" \cite{Hallard_1991}. 
Smith et al. addressed this problem by presenting solutions for concave diagram tiles without 
matching conditions \cite{Smith_2024a, Smith_2024b}. According to Smith et al. \cite{Smith_2024a}, 
\mbox{Tile$(a, b)$}, which can generate only non-periodic tiling, is an ``aperiodic monotile," 
in which all tiles are of the same size and shape; however, the tiling generation 
process requires the use of reflected tiles. According to Smith et al. \cite{Smith_2024b}, 
the tiles ``Spectres," which can generate only non-periodic tiling, do not require the 
use of reflected tiles during the tiling generation process. 
The tiles are referred to as ``chiral aperiodic monotiles." \footnote{ 
The monotile \mbox{Tile$(a, b)$} with specific values of $a$ and $b$ can generate periodic tilings. 
For example, \mbox{Tile$(1, 1)$}, corresponding to $a = b = 1$ in \mbox{Tile$(a, b)$}, admits  
periodic tilings. By contrast, \mbox{Tile$(1, 1)$} can generate only non-periodic tilings if and 
only if it does not allow the use of reflected tiles during the tiling generation process, and is 
referred to as a ``weakly chiral aperiodic monotile." In comparison, Spectres can generate 
only non-periodic tilings using only one side, either anterior or posterior, even if they allow 
the use of reflected tiles during the tiling generation process \cite{Smith_2024a, Smith_2024b}.
}

Using the fact that every triangle and quadrilateral is a monotile, that convex hexagonal 
monotiles can be classified into three Type families each of which admits a periodic tiling, 
and that no convex polygonal monotiles with seven or more edges exist  \cite{Bollo_1963, 
Gardner_1975a, G_and_S_1987, Niven_1978, Klam_1980, Reinhardt_1918, Zong_2020},
together with Theorem~\ref{thm1}, we obtained the following result \cite{Sugimoto_2012, 
Sugimoto_2015a, Sugimoto_2015b, Sugimoto_2016, Sugimoto_2017a}.

\begin{thm}\label{thm2}
In the convex polygonal monotiles that can generate an edge-to-edge tiling, no aperiodic 
monotile exists \rm{(}\textit{without matching conditions other than edge-to-edge}\rm{)}. 
\rm{(}\textit{In other words, without matching conditions other than edge-to-edge, 
no convex polygonal monotile can be an aperiodic monotile.}\rm{)}
\end{thm}

In \cite{Rao_2017}, Rao claimed that there is no convex polygon that allows only 
non-periodic tilings (without matching conditions).

\section{}
\label{appD}

Let $P$ be a tile on the Euclidean plane. A \textit{corona} of $P$ is the collection 
consisting of the centrally placed original $P$ and its surrounding layer of congruent 
copies of $P$. The corona is formed without gaps or overlaps (except at the 
boundaries of $P)$. The first corona is the collection of all tiles that share 
a boundary point with a centrally placed tile (including the original tile 
itself) \cite{Goodman_2017, G_and_S_1987, Sugimoto_2018a}.

\begin{figure}[htbp]
\centering\includegraphics[width=15cm,clip]{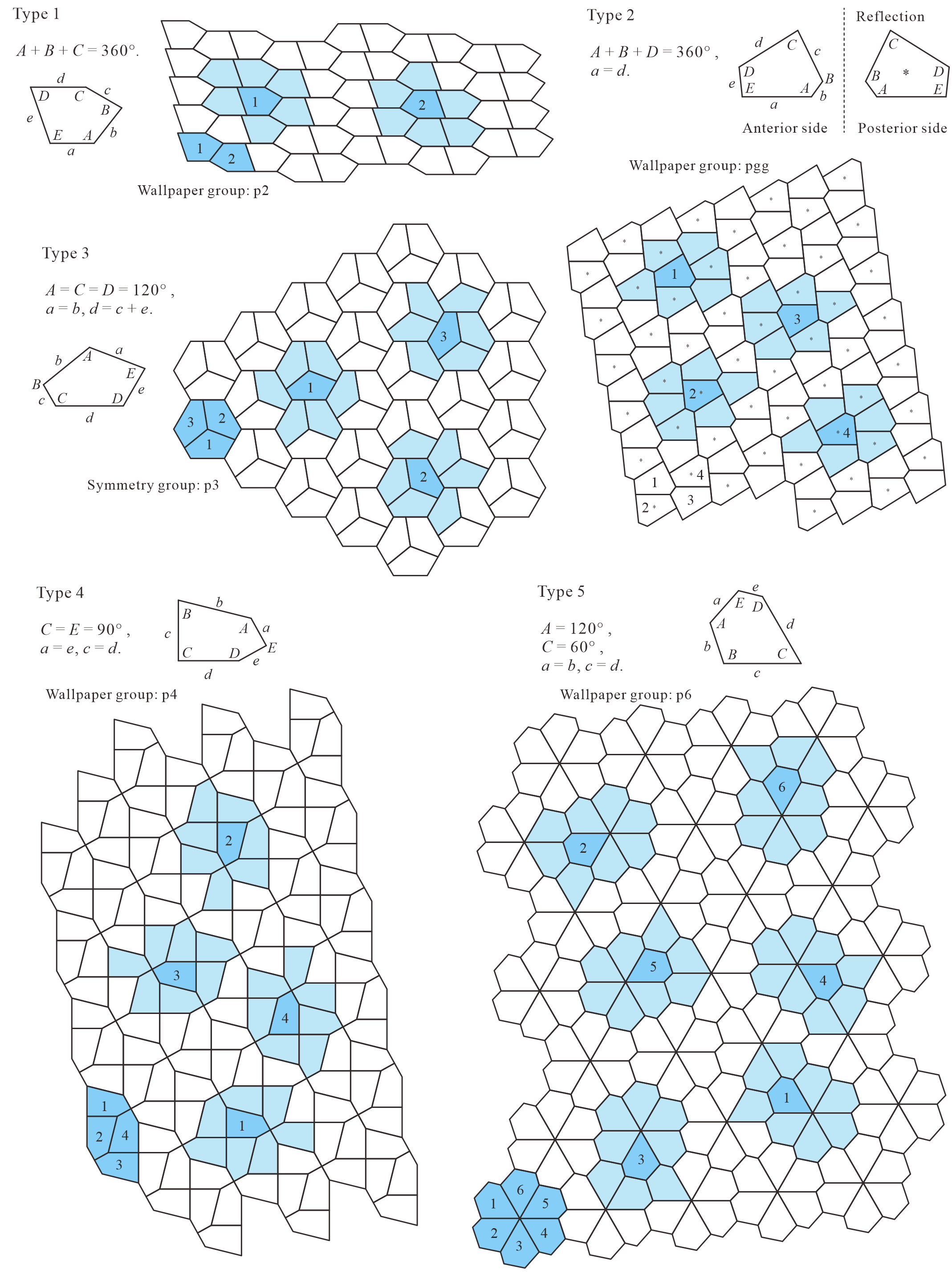} 
\caption{{\small 
Isohedral representative tilings of Types 1--5 and the first coronas of 
each convex pentagonal monotile that constitutes their translation units. 
In this figure, the convex pentagons corresponding to the reflected 
tiles are marked with an asterisk ``*.''
}
\label{fig13}
}
\end{figure}

\begin{figure}[htbp]
\centering\includegraphics[width=15cm,clip]{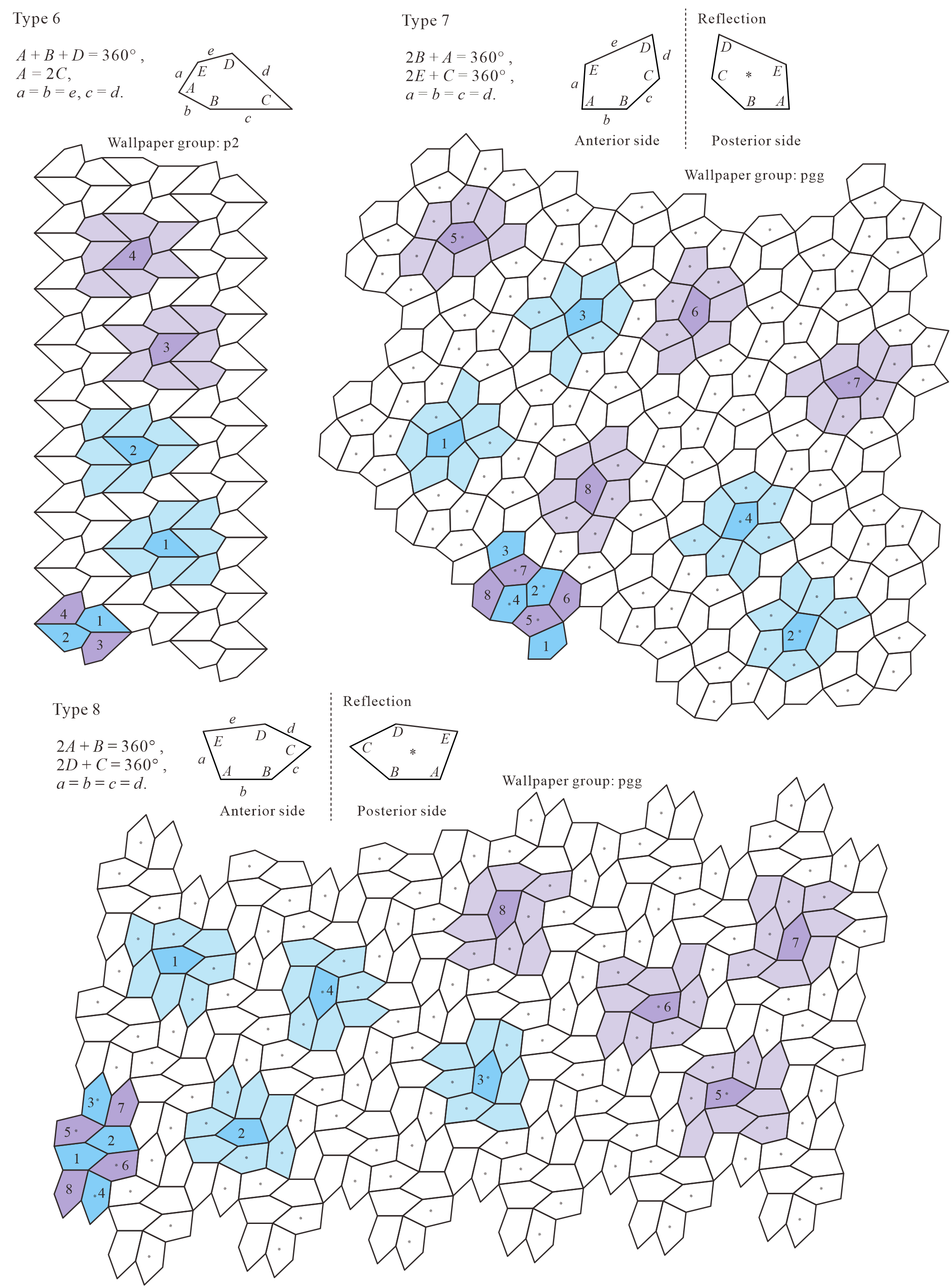} 
\caption{{\small 
2-isohedral representative tilings of Types 6--8 and the first coronas 
of each convex pentagonal monotile that constitutes their translation units. 
In this figure, the convex pentagons corresponding to the reflected 
tiles are marked with an asterisk ``*.''
}
\label{fig14}
}
\end{figure}

\begin{figure}[htbp]
\centering\includegraphics[width=15cm,clip]{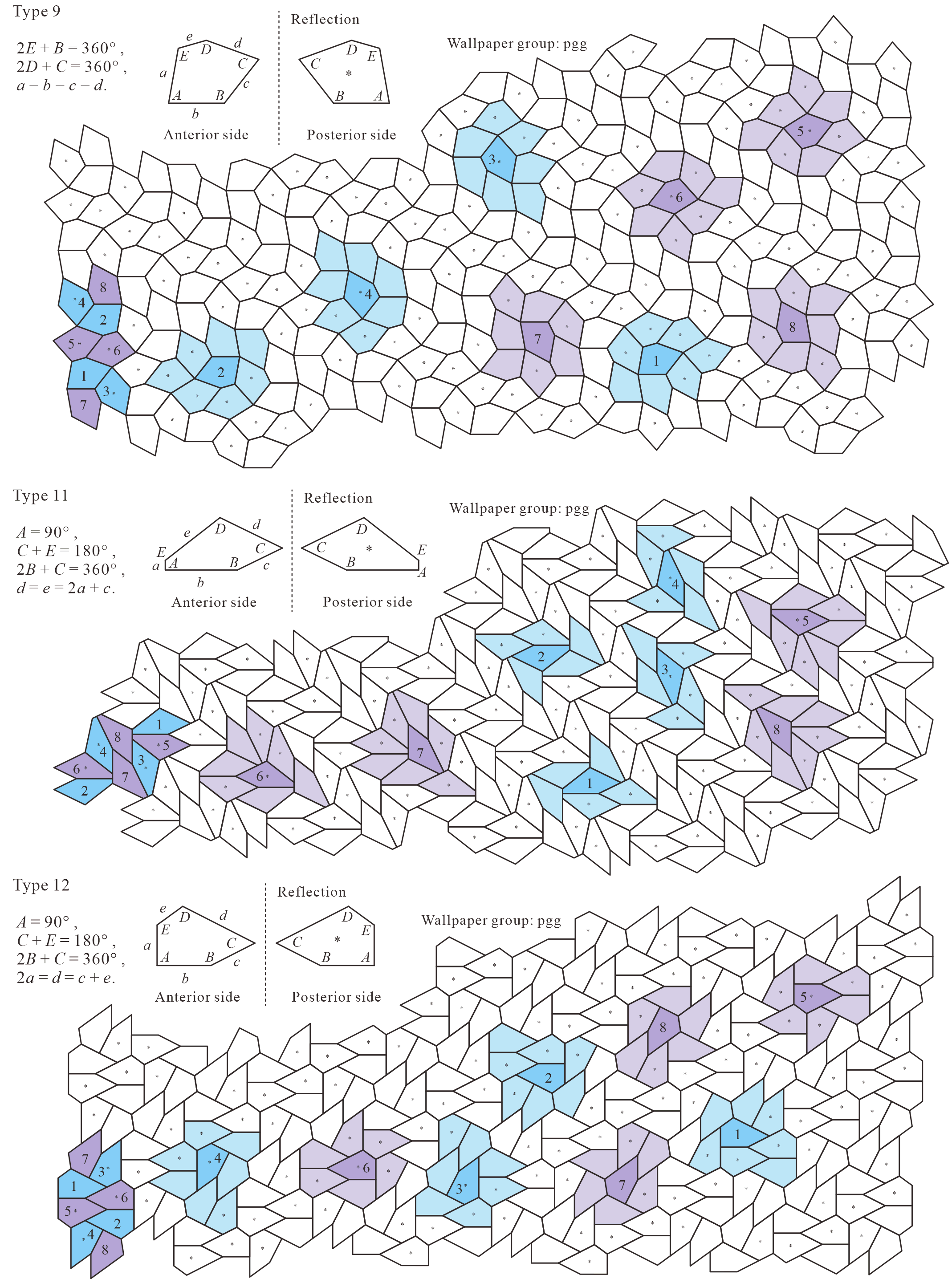} 
\caption{{\small 
2-isohedral representative tilings of Types 9, 11, and 12 and the first coronas 
of each convex pentagonal monotile that constitutes their translation units. 
In this figure, the convex pentagons corresponding to the reflected 
tiles are marked with an asterisk ``*.''
}
\label{fig15}
}
\end{figure}

\begin{figure}[htbp]
\centering\includegraphics[width=15cm,clip]{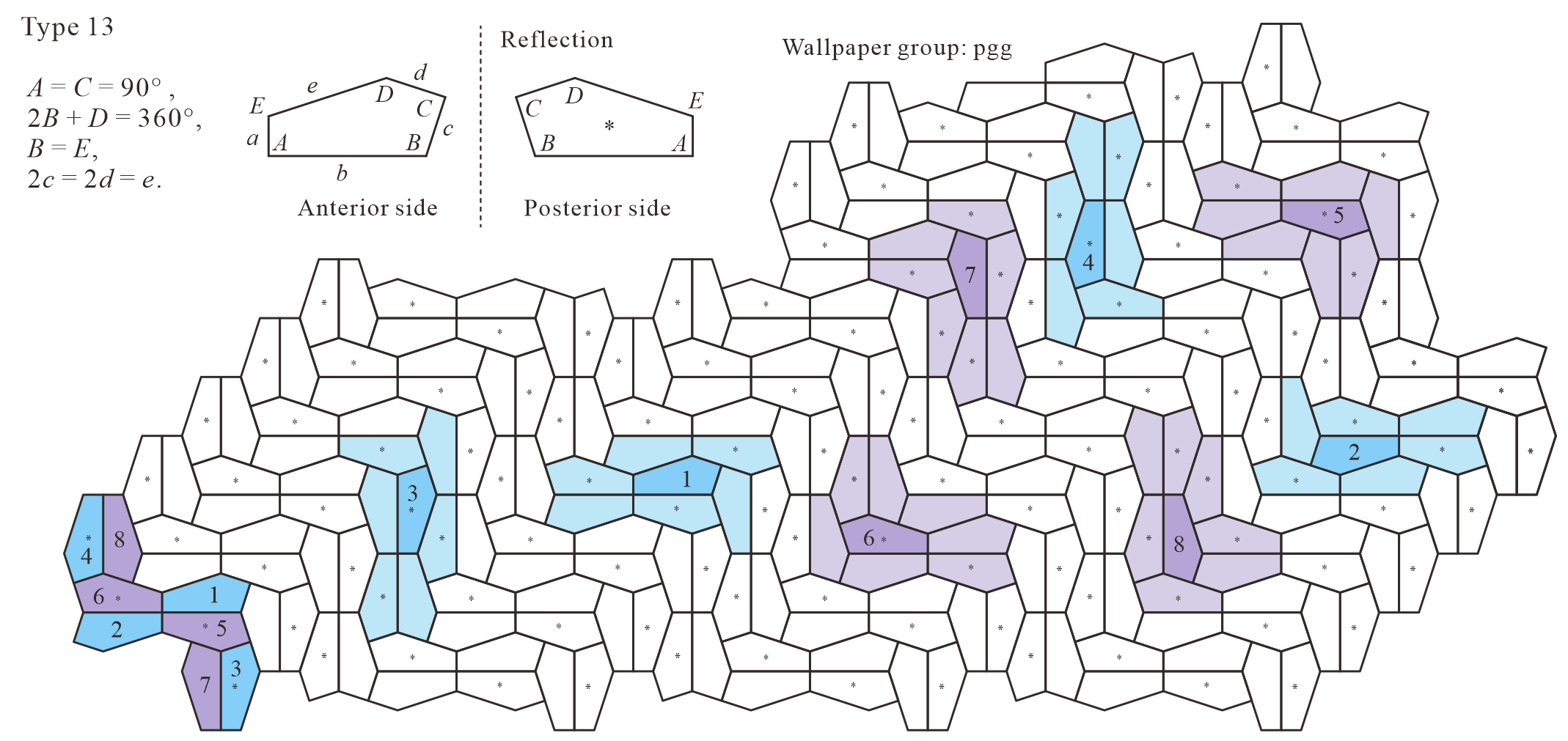} 
\caption{{\small 
2-isohedral representative tiling of Types 13 and the first coronas 
of each convex pentagonal monotile that constitutes the translation unit. 
In this figure, the convex pentagons corresponding to the reflected 
tiles are marked with an asterisk ``*.''
}
\label{fig16}
}
\end{figure}

\begin{figure}[htbp]
\centering\includegraphics[width=15cm,clip]{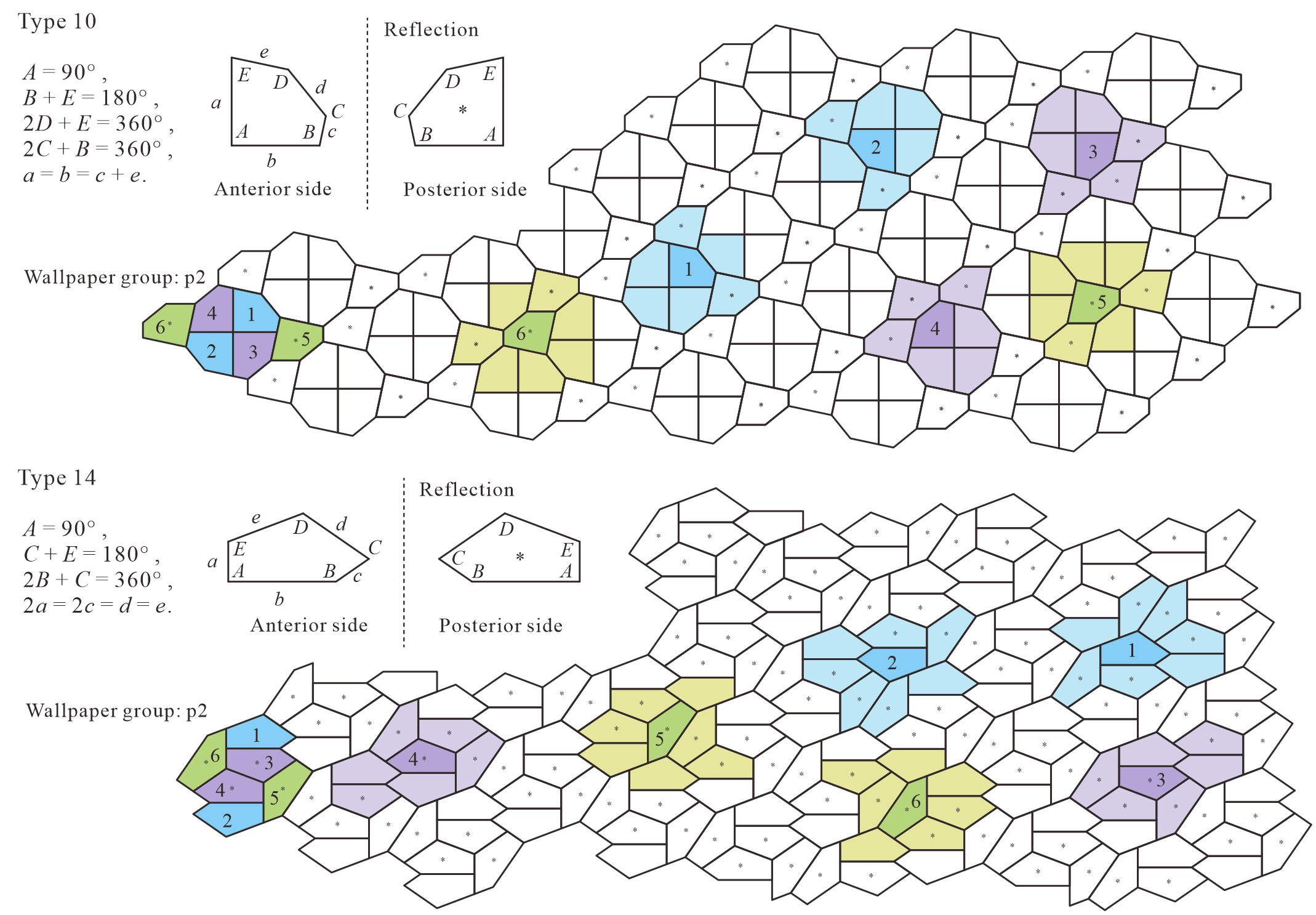} 
\caption{{\small 
3-isohedral representative tilings of Types 10 and 14 and the first coronas 
of each convex pentagonal monotile that constitutes their translation units. 
In this figure, the convex pentagons corresponding to the reflected 
tiles are marked with an asterisk ``*.''
}
\label{fig17}
}
\end{figure}

\begin{figure}[htbp]
\centering\includegraphics[width=15cm,clip]{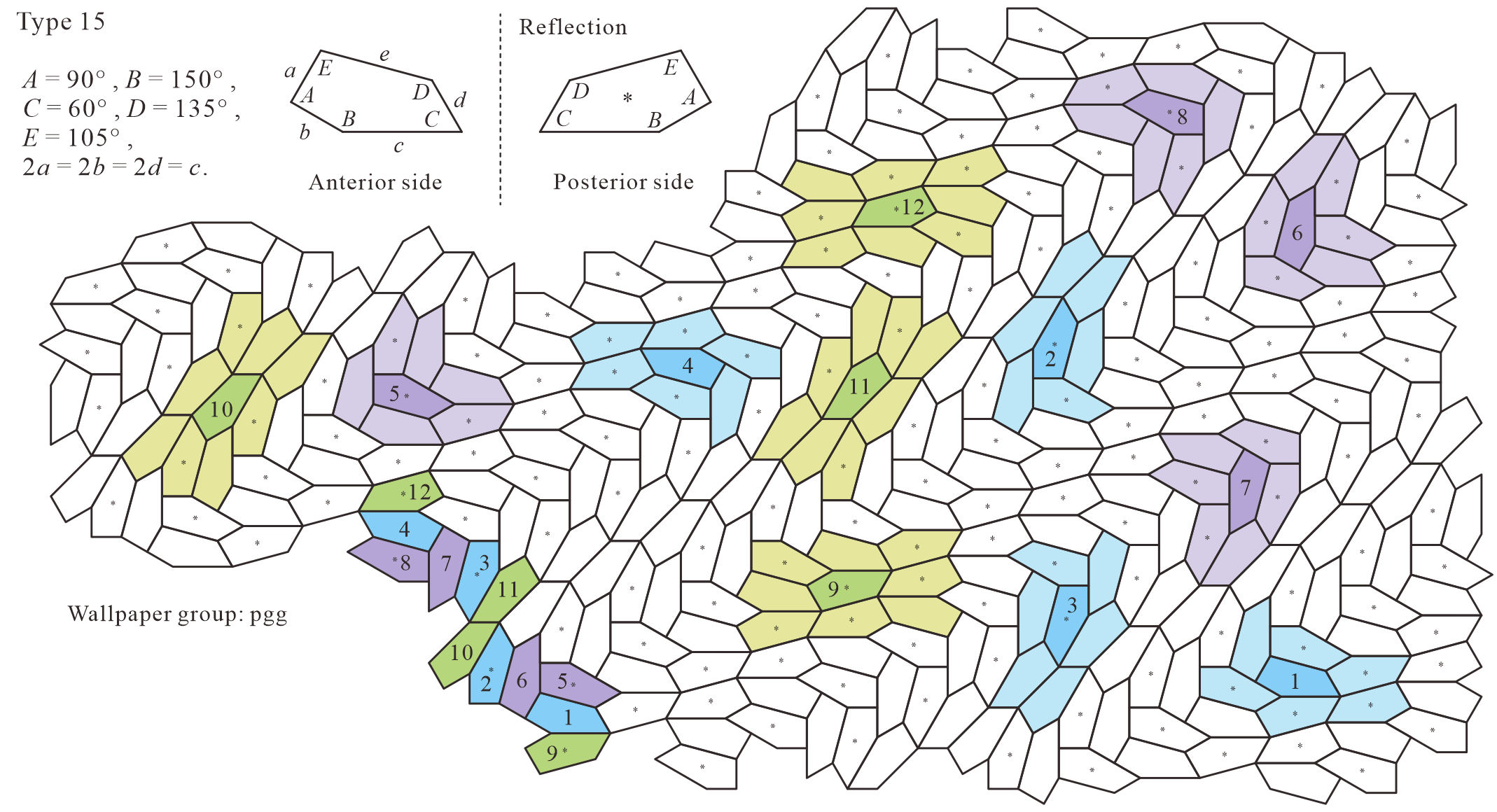} 
\caption{{\small 
3-isohedral representative tiling of Type 15 and the first coronas of 
each convex pentagonal monotile that constitutes the translation unit. 
In this figure, the convex pentagons corresponding to the reflected 
tiles are marked with an asterisk ``*.''
}
\label{fig18}
}
\end{figure}

The first coronas of all tiles in isohedral tilings have the same shape. 
(Remark: If coronas have the same shape as their reflected counterparts, 
they are considered to be of the same type.) Therefore, in the representative 
tilings of Types 1--5 shown in Figure~\ref{fig01}, which correspond to 
isohedral tilings, the shapes of the first coronas of each convex pentagonal 
monotile forming the translation units are of a single type (see Figure~\ref{fig13}). 
In contrast, in the representative tilings of Types 6--9 and 11--13 shown in 
Figure~\ref{fig01}, which correspond to 2-isohedral tilings, the shapes of the first 
coronas of each convex pentagonal monotile forming the translation units are 
of two distinct types (see Figures~\ref{fig14}--\ref{fig16}). Furthermore, in the 
representative tilings of Types 10, 14, and 15 shown in Figure~\ref{fig01}, which 
correspond to 3-isohedral tilings, the shapes of the first coronas of each 
convex pentagonal monotile forming the translation units are of three 
distinct types (see Figures~\ref{fig17} and \ref{fig18}).
The symbols such as ``p2," labeled as ``Wallpaper group: p2" and appearing 
near the representative tilings of each Type in Figures~\ref{fig13}--\ref{fig18}, 
denote the wallpaper groups (i.e., plane symmetry groups or plane crystallographic 
groups) of the depicted tilings\footnote{ 
In certain special cases of convex pentagons, the wallpaper groups of the 
representative tilings may change. For example, when the representative tiling of 
Type 10 is formed using the convex pentagonal monotile in $T_{1} \cap T_{10}$ 
as shown in Figure~\ref{fig04}, the resulting tiling has the wallpaper group cmm 
\cite{G_and_S_1987, Schattschneider_1978, wiki_pentagon_tiling}.
} 
\cite{G_and_S_1987, Schattschneider_1978, wiki_pentagon_tiling, wiki_wallpaper_group}.

\begin{figure}[htbp]
\centering\includegraphics[width=15cm,clip]{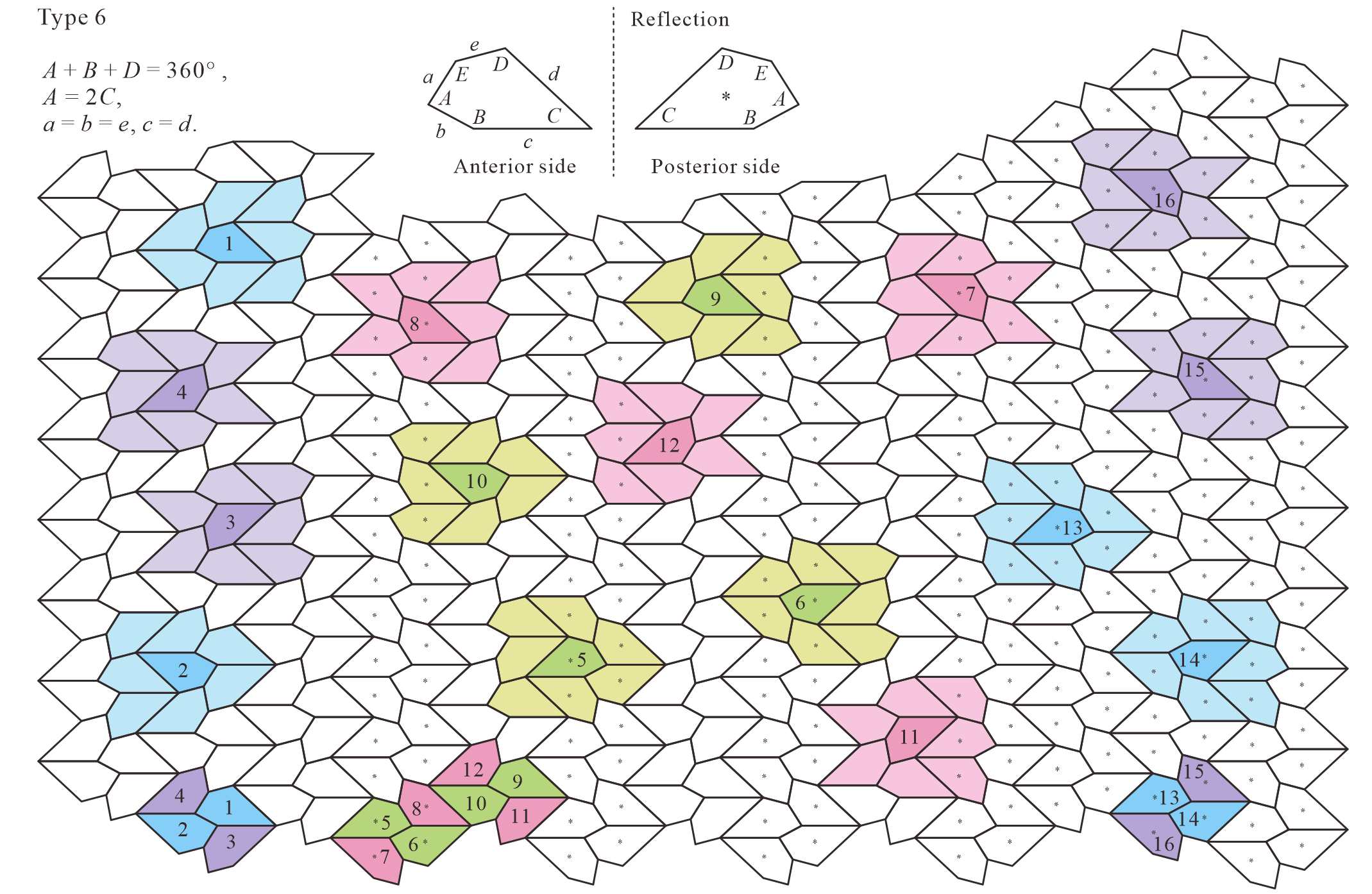} 
\caption{{\small 
4-isohedral tiling generated by convex pentagonal monotiles 
belonging only to the Type 6 family. In this figure, the convex pentagons 
corresponding to the reflected tiles are marked with an asterisk ``*.''
}
\label{fig19}
}
\end{figure}

All convex pentagonal monotiles in $T_{6}$ can generate edge-to-edge tilings, 
as shown in Figure~\ref{fig19}, using reflected tiles, because the belts they form 
through vertical translation in the same direction can be freely connected 
horizontally. As shown in Figure~\ref{fig19}, convex pentagonal monotiles in 
$T_{6} \cap T_{1}^{C}$ (i.e., convex pentagonal monotiles belonging to 
the Type 6 family that are not contained in $T_{1})$ can generate 4-isohedral 
tilings owing to the free combinatorial properties of these belts when 
reflected tiles are used\footnote{ 
Based on the Venn diagram in Figure~\ref{fig03}, a convex pentagonal monotile in 
$T_{1} \cap T_{6}$ corresponds to the convex pentagon in $T_{1} \cap T_{5} \cap T_{6}$ 
as shown in Figure~\ref{fig04}. The intersection ``$T_{1} \cap T_{5} \cap T_{6}$'' 
contains one fixed shape that exhibits line symmetry. Due to this property, the 
shapes of the first coronas are still of two types (i.e., the tiling is 2-isohedral) 
when it generates the tiling shown in Figure~\ref{fig19} using the reflected tile. 
Therefore, the case of the pentagon in $T_{1} \cap T_{6}$ (i.e., $T_{1} \cap T_{5} \cap T_{6}$) 
is excluded.
}. 
Based on the above properties and the Venn diagram in Figure~\ref{fig03}, and 
considering that convex pentagonal monotiles that can generate isohedral tilings 
belong to at least one of the Type 1--5 families shown in Figure~\ref{fig01} \cite{G_and_S_1987, 
Mann_2015, Scherphuis, wiki_pentagon_tiling, Zong_2020}, convex pentagonal monotiles 
that belong only to the Type 6 family (i.e., convex pentagonal monotiles in 
$\bigl( T_{6} \cap T_{1}^{C} \bigr) \cap T_{2}^{C})$ \footnote{ 
According to the Venn diagram in Figure~\ref{fig03}, the correct expression should be 
``convex pentagonal monotile in $\bigl( \bigl( T_{6} \cap T_{1}^{C} \bigr) \cap T_{2}^{C} \bigr) \cap T_{5}^{C}$." 
However, the convex pentagonal monotile in $T_{5} \cap T_{6}$ corresponds solely 
to the convex pentagon in $T_{1} \cap T_{5} \cap T_{6}$, as shown in Figure~\ref{fig04}. 
Therefore, ``convex pentagonal monotile in $\bigl( T_{6} \cap T_{1}^{C} \bigr) \cap T_{2}^{C}$" and 
``convex pentagonal monotile in $\bigl( \bigl( T_{6} \cap T_{1}^{C} \bigr) \cap T_{2}^{C} \bigr) \cap T_{5}^{C}$" 
are identical.
} 
are anisohedral tiles that can generate two types of tilings: 2-isohedral and 
4-isohedral. Furthermore, they are classified as convex polygonal 
anisohedral tiles that can generate multiple types of $k$-isohedral tilings.

It is possible that convex polygonal anisohedral tiles that can generate 
multiple types of $k$-isohedral tilings are exclusively convex pentagonal 
monotiles that belong only to the Type 6 family. First, there are no anisohedral 
tiles for triangular, convex quadrilateral, or convex hexagonal monotiles, and no 
convex polygonal monotiles with seven or more edges \cite{Bollo_1963, 
G_and_S_1987, Niven_1978, Klam_1980, Reinhardt_1918, Zong_2020}. 
Next, if Rao's claim \cite{Rao_2017} is correct, there are no convex polygonal 
anisohedral tiles other than convex pentagonal monotiles that are not contained in 
$T_{1}$, $T_{2}$, $T_{3}$, $T_{4}$, or $T_{5}$. Consequently, considering the Venn 
diagram in Figure~\ref{fig03} and the list in Figure~\ref{fig04}, it follows that there are 
no convex polygonal anisohedral tiles other than convex pentagonal monotiles 
belonging to the Types 6--15 families that are not contained in $T_{1}$ or 
$T_{2}$. Furthermore, convex pentagonal monotiles that belong to the Types 
7--15 families require reflected tiles to form the representative tilings of 
each Type. (Unlike the convex pentagonal monotiles that belong to the Type 6 
family, the representative tilings use reflected tiles.) Based on these 
observations, we conjecture that it is unlikely that convex pentagonal 
monotiles belonging to the Types 7--13 families, which are not contained in 
$T_{1}$ or $T_{2}$, can generate multiple types of $k$-isohedral tilings. 

To verify whether this conjecture is correct, it is necessary to determine whether all 
convex pentagonal monotiles in $\bigl( T_{w} \cap T_{1}^{C} \bigr) \cap T_{2}^{C}$ 
(i.e., convex pentagonal monotiles belonging only to the Type $w$ family\footnote{ 
Considering the Venn diagram in Figure~\ref{fig03} and the list in Figure~\ref{fig04}, 
convex pentagonal monotiles in $\bigl( T_{w} \cap T_{1}^{C} \bigr) \cap T_{2}^{C}$ 
(i.e., convex pentagonal monotiles belonging to the Type $w$ family, which are not 
contained in $T_{1}$ or $T_{2})$ are convex pentagonal monotiles belonging only 
to the Type $w$ 
family.}, 
where $w = 6-13$)  can generate tiling patterns other than the representative 
tilings of Type $w$. If such tiling patterns exist, further investigation into their 
properties will be required (note that this investigation is similar to that 
described in \ref{appB}). We believe that investigating convex pentagonal 
monotiles with $w = 6$ (i.e., those that belong only to the Type 6 family) 
is important, as there may be cases analogous to the convex pentagonal 
monotile in $T_{1} \cap T_{5} \cap T_{6}$.

Some examples of 2-isohedral, 3-isohedral, 4-isohedral, and 6-isohedral 
tilings formed by convex pentagonal monotiles are presented in \cite{Scherphuis}. 
However, the convex pentagonal monotiles that form these tilings appear to 
be able to generate isohedral tilings as well (i.e., they are not 
anisohedral tiles).

\renewcommand{\figurename}{{\small Figure}}

\renewcommand{\figurename}{{\small Figure}}

\renewcommand{\figurename}{{\small Figure}}

\renewcommand{\figurename}{{\small Figure}}

\renewcommand{\figurename}{{\small Figure}}

\renewcommand{\figurename}{{\small Figure}}

\renewcommand{\figurename}{{\small Figure}}

\end{document}